\documentclass[11pt,reqno]{amsart}
 \usepackage{amsgen, amstext,amsbsy,amsopn, amsthm, amsfonts,amssymb,amscd,amsmath,euscript,enumerate,url,verbatim,calc,xypic,tikz}
\usepackage{hyperref}
\usepackage{MnSymbol}
\usepackage{pgfkeys}
\usepackage{mathtools}
\usepackage{eqnarray,amsmath}

\let\emptyset\varnothing

\def\multiset#1#2{\ensuremath{\left(\kern-.3em\left(\genfrac{}{}{0pt}{}{#1}{#2}\right)\kern-.3em\right)}}

\usetikzlibrary{arrows}
\newcommand{\midarrow}{\tikz \draw[-triangle 90] (0,0) -- +(.1,0);}
 \usepackage{latexsym}
 \usepackage{graphics}
 \usepackage{color}
\usepackage{lastpage}
\usepackage{fancyhdr}
\usepackage{multirow}
\allowdisplaybreaks
\usepackage{graphicx}
\graphicspath{ {F:/IMAGES/} }

\newcommand{\PP}{\mathcal{P}}

  \newcommand{\pd}{\operatorname{pd}}

  \newcommand{\Ass}{\operatorname{Ass}}

 \newcommand{\Supp}{\operatorname{Supp}} 
\newcommand{\rad}{\operatorname{rad}}
  \newcommand{\reg}{\operatorname{reg}}

\newcommand{\proset}{\,\mathrel{\lower 4pt\hbox{$\scriptscriptstyle/$}
\mkern -14mu\subseteq }\,} 

 \newtheorem{theorem}{Theorem}[section]
 \newtheorem{corollary}[theorem]{Corollary}
 \newtheorem{lemma}[theorem]{Lemma}
 \newtheorem{proposition}[theorem]{Proposition}

 \newtheorem{question}{Question}

\usepackage{amsmath}
\newtheorem{notation}[theorem]{Notation}

 \theoremstyle{definition}
 
 \newtheorem{remark}[theorem]{Remark}
 \newtheorem{definition}[theorem]{Definition}

 \newtheorem{example}[theorem]{Example}
\title[Monomial Ideals of VEW graphs] {Monomial Ideals of VEW graphs}
\author[ M. Kumar and D.K. Pradhan ]{Manohar Kumar and Dipak Kumar Pradhan$^*$  }

\thanks{AMS Classification 2010: 05C22, 13D02, 13D45, 05C25, 05C38, 05E40}
\thanks{$^*$ Corresponding author}
\address{Department of Mathematics, Indian Institute of Technology
		Madras, Chennai, Tamil Nadu, INDIA - 600036.}
\email{manhar349@gmail.com}

\address{Department of Mathematics, S.N. Degree College, Kendrapara, Utkal University,  Odisha, INDIA - 754220.}\email{dipakmaths123@gmail.com}

\begin{document}
\maketitle

\begin{abstract}
Let $D$ be a $k$-edge weighted VEW graph and $I=I(D)$ be its edge ideal. In this paper,  we characterize the irredundant irreducible decomposition of edge ideal $I$. Also, we study the Cohen-Macaulay property and Castelnuovo-Mumford regularity of different classes of $k$-edge weighted VEW graphs.    
\end{abstract}

\section{Introduction}
Monomial ideals provide a link among Commutative Algebra, Combinatorics, and Topology. The squarefree monomial ideals have been extensively studied through {\em Stanley-Reisner ideals} of {\em simplicial complexes} and {\em edge ideals} of {\em hypergraphs}. In particular, quadratic squarefree monomial ideals are studied via edge ideals of simple graphs. Recently, much research has focused on non-squarefree monomial ideals through edge ideals of {\em vertex-weighted oriented graphs} and {\em edge-weighted graphs}. See \cite{kanoy1, kblo22, cs13,cpr22, fsty21, gmsv18,  kumar1, kumar2, kumar3, mandal1, mandal2, mandal3, mandal4,pitones, zhu}
 for details. In this paper, we introduce a more general class of non-squarefree monomial ideals via {\em vertex-edge weighted oriented graphs}.
\vspace{0.2cm}\\
A \textit{directed graph} or\textit{ digraph} $D=(V(D),E(D))$ consists of a finite set $V(D)$ of vertices together with a prescribed collection $E(D)$ of ordered pairs of distinct vertices called edges or arrows. If $ (u,v) \in E(D),$ then we call it a directed edge, where the direction is from $ u $ to $ v $ and $ u $ (respectively $ v $) is called the initial vertex (respectively the terminal vertex). An oriented graph is a directed graph without multiple edges or loops. In other words, an oriented graph $D$ is a simple graph $G$ together with an orientation of its edges. 
An oriented graph $D$ is called vertex-edge weighted oriented if it is equipped with two weight functions $w:V(D)\longrightarrow \mathbb N$ and $w':E(D) \longrightarrow \mathbb N $. In fact, a vertex-edge weighted oriented graph (notationally VEW graph) $D$ is a quadruplet $D=(V(D), E(D), w,w')$, where  $V(D)$ and  $E(D)$ are the vertex set and edge set,  respectively, the weight of a vertex $x_i\in V(D)$ is $w(x_i)$, denoted by $w_i$ and the weight of an edge $(x_i,x_j) \in E(D)$ is $w((x_i,x_j))$ denoted by $w_{(x_i,x_j)}$. If $V(D)=\{x_1,\ldots ,x_n\},$ then we can regard each vertex $x_i$ as a variable and consider the polynomial ring $R=\mathbb{K}[x_1,\ldots ,x_n]$ over a field $\mathbb{K}$. 
 Then the {\em edge ideal } of $D$  is defined  as
$$I(D)  =  (x_i^{\min\{w_{(x_i,x_j)},w(x_i)\}} x_j^{\max\{w_{(x_i,x_j)},w(x_j)\}}~|~(x_i,x_j) \in E(D) ).$$ The underlying graph of $ D $ is the simple graph $ G $ whose
vertex set is $V(G) = V(D) $ and whose edge set is $E(G)= \left\{
\{u,v\}~|~(u, v) \in E(D)  \right\}.$  The edge ideal of $ G $ is 
$I(G) = ( uv ~|~  \{u, v\} \in E(G) ).$
If $w_{(x_i,x_j)} = 1 $ for all $(x_i, x_j) \in E(D)  $ and $w(x_i)=1$ for all $x_i \in V(D)$, then $I(D) = I(G)$.    
 Setting all edge weights to $1$, VEW edge ideal yields the vertex weighted oriented ideal, which was initially studied in \cite{pitones}.

In this paper, we study {\rm $k$-edge weighted} VEW graphs, that is, VEW
graphs with all edges of weight $k$. Then the edge ideal of a $k$-edge weighted VEW graph $D$ is defined as
$$ I(D)  =  (x_i^{\min\{k,w(x_i)\}} x_j^{\max\{k,w(x_j)\}}~|~(x_i,x_j) \in E(D) ).$$

If a vertex $x_i$ of $D$ is a source (i.e., it has only arrows leaving $x_i$), then we always assume $w(x_i)\leq k$. Indeed, if $w(x_i)>k$, then for every edge $(x_i,x_j)\in E(D)$ we have 
$x_i^{\min\{k,w(x_i)\}}x_j^{\max\{k,w(x_j)\}}=x_i^{k}x_j^{\max\{k,w(x_j)\}}$, because $\min\{k,w(x_i)\}=k$. Hence the generators of $I(D)$ involving $x_i$ do not depend on the exact value of $w(x_i)$ when $w(x_i)>k$. Therefore, we always assume $w(x_i)\le k$. Similarly, if a vertex $x_i$ of $D$ is a sink (i.e., it has only arrows coming towards $x_i$), we always assume $w(x_i)\ge k$. Indeed, if $w(x_i)<k$, then for every edge $(x_j,x_i)\in E(D)$ we have
$x_j^{\min\{k,w(x_j)\}}x_i^{\max\{k,w(x_i)\}}=x_j^{\min\{k,w(x_j)\}}x_i^{k}$, since $\max\{k,w(x_i)\}=k$. Hence the generators of $I(D)$ involving $x_i$ do not depend on the exact value of $w(x_i)$ when $w(x_i)<k$. Therefore, we always assume $w(x_i)\ge k$.  If $k=1$, then $ I(D)  =  (x_i x_j^{w(x_j)}~|~(x_i,x_j) \in E(D))$ is called edge ideal of $1-$edge weighted VEW graph. For edge ideal of any $1$-edge weighted VEW graph (i.e., weighted oriented graph), unmixedness, Cohen-Macaulayness, componentwise linearity, Castelnuovo-Mumford regularity and symbolic powers have been studied extensively in \cite{kanoy1, kblo22, cpr22, gmsv18,   kumar1,kumar2, kumar3, mandal1,mandal2,mandal3,mandal4, pitones, zhu}.	
\vspace{0.2cm}\\
In this paper, we provide irreducibe decomposition of edge ideal of a $k$-edge weighted VEW graph (Theorem \ref{irreducible}). Furthermore, we give criteria for the unmixedness and Cohen-Macaulayness of the edge ideal of certain classes of $k$-edge weighted VEW graphs (Theorem \ref{cm-forest}, Theorem \ref{cm-graph}), by extending the results from \cite[Theorem 3.1]{H.T-2018} and \cite[Theorem 5]{gmsv18} to any $k$-edge weighted VEW graph. Finally, we provide formulas for the regularity of certain classes of $k$-edge weighted VEW graphs (Theorem \ref{atmost-one}, Theorem \ref{ClassII}, Theorem \ref{cycle1}, Theorem \ref{cycle2}).
\vspace{0.2cm}\\
Now, we give a section-wise description.
In Section \ref{sec2}, we recall the definitions and the basic results that are required to prove our main results. In Section \ref{sec3}, we discuss strong vertex cover of $k$-edge weighted VEW graphs. In  Section \ref{sec4}, we identify the irredundant irreducible decomposition of any $k$-edge  weighted VEW graph $D$. We show that the minimal monomial irreducible ideals of $I(D)$ are associated with the strong vertex covers of $D$. In  Section \ref{sec5}, we study Cohen-Macaulayess of $k$-edge weighted VEW forest and Cohen-Macaulaynes of $k$-edge  weighted VEW graph with perfect matching $\{\{x_1,y_1\},\ldots,\{x_r,y_r\}\}$. Finally, in Section \ref{sec6}, using the technique of labeled hypergraph, we give one formula to compute regularity of naturally oriented VEW paths, cycles, rooted forests, unicyclic graphs, etc.  Also, we find regularity of $k$-edge weighted VEW paths and cycles assuming non-source vertices with weight $\neq k  $  are sinks. In the last Section \ref{sec7}, we discuss some future problems.

\section{Preliminaries} \label{sec2}           
In this section, we recall some algebraic definitions and results. Let $\mathbb{K}$ be a field and $R=\mathbb{K}[x_1,\ldots ,x_n]$ be the polynomial ring in $n$ variables over $\mathbb{K}$. Suppose that $M$ is a non zero graded $R$-module with finite graded minimal free resolution $$0 \longrightarrow  \cdots   \longrightarrow \underset{j} \bigoplus  R(-j)^{\beta_{1,j}(M)} \longrightarrow \underset{j} \bigoplus  R(-j)^{\beta_{0,j}(M)}  \longrightarrow    M    \longrightarrow    0,$$
where  $ {\beta_{i,j}}(M) $ denote  the ($i,j$)-th graded Betti number of $ M$. The invariant which measures the complexity of the module is the Castelnuovo-Mumford regularity denoted by $\reg(M)$ and  defined as
\begin{equation*}
	\reg(M) := \max \{ j -  i \mid  {\beta_{i,j}}(M)  \neq  0 \} .
\end{equation*}

 With the help of Betti splitting, we can compute this type of invariant for certain classes of ideals. The Betti splitting is defined in the following.

\begin{notation}
  Let $I \subset R$ be a monomial ideal. Then $\mathcal{G}(I)$ denotes the set of minimal monomial generators of $ I $.  
\end{notation}
 
\begin{definition}\label{2.1}  
	Let $I$ be a monomial ideal and suppose that there exist monomial ideals $J$ and $K$ such that $\mathcal{G}(I)$ is the disjoint union of $\mathcal{G}(J)$   and $\mathcal{G}(K)$. Then
	$I = J + K$ is a Betti splitting if
	\begin{equation*}
		{\beta_{i,j}}(I) =  {\beta_{i,j}}(J) +
		{\beta_{i,j}}(K) +  {\beta_{i-1,j}}(J\cap K)
	\end{equation*}
	for all $i, j  \geq  0,$ 	where $ {\beta_{i-1,j}}(J\cap K) = 0 $ if $ i = 0.$
\end{definition}
This formula was first obtained for the total Betti numbers by Eliahou and Kervaire \cite{elia}. In \cite{francisco}, the authors describe the following sufficient conditions for an ideal $I$ to have a Betti splitting.
\begin{theorem}\cite[Corollary 2.7]{francisco} \label{betti.1}	Suppose that $I = J + K$, where $\mathcal{G}(J)$ contains
	all the generators of $I$ divisible by some variable $x_i$ and $\mathcal{G}(K)$ is a nonempty set
	containing the remaining generators of $I$. If $J$ has a linear resolution, then $I = J+K$
	is a Betti splitting.
\end{theorem}
When $I$ has a Betti splitting, Definition \ref{2.1} implies the following result:
\begin{corollary}\label{betti.2}
	If $I = J + K$ is a Betti splitting, then
	$$ \reg (I)  = \max \{ \reg (J),~ \reg (K),~ \reg (J  \cap  K)  -  1 \} .$$
\end{corollary}

Let $u \in R$ be a monomial, we set $\Supp(u) = \{x_i : x_i  ~|~ u \}.$ Let $I$ be a monomial ideal,
$\mathcal{G}(I) = \{u_1,\ldots, u_m\}$ denote the unique minimal set of monomial generators of $I$ and we set $\Supp(I) :=  \displaystyle{\bigcup_{i=1}^m}
\Supp(u_i).$ Now, we give some useful lemmas on regularity.      

\begin{lemma}\cite[Lemma 3.2]{ht10}\label{reg.2}
	Let $R_1 = \mathbb{K}[x_1,\ldots , x_m]$ and $R_2 = \mathbb{K}[x_{m+1}, \ldots, x_n]$ be
	two polynomial rings,  $ I \subset $ $R_1$ and  $ J \subset $ $R_2$ be two nonzero homogeneous ideals. Then
	\begin{enumerate}
		\item $ \reg (I + J) = \reg (I) + \reg (J) - 1.$
        \item $ \reg (IJ) = \reg (I) + \reg (J).$
	\end{enumerate}
\end{lemma}

\begin{lemma}\cite[ Lemma 1.2]{H.T-2016}  \label{reg.4}  Let  $0  \rightarrow A \rightarrow B \rightarrow C \rightarrow 0$  be 	short exact sequence of finitely generated graded $R$-modules. Then
	\item $\reg(B)$ $\leq \max \{\reg (A),~ \reg (C)\}$ and the equality holds if $\reg(A)$ $- 1  \neq$ $\reg(C).$ 
\end{lemma}
\begin{lemma}  \cite[Lemma 3.1]{htha}\label{H<G}
	Let $G = (V, E)$ be a simple graph. If $ G^\prime $ is an induced subgraph of $ G $, then $ \reg(I(G^\prime)) \leq  \reg(I(G)). $	
\end{lemma}
The propositions below are based on the regularity of the edge ideals of a path and a cycle. Recall that a \textit{path} of length $n$ is denoted by $P_n$, and a \textit{cycle} of length $n$ (or $n$-\textit{cycle}), denoted by $C_n$, is a connected graph on the vertex set $V(C_n)=\{x_1,\ldots,x_n\}$ such that $E(C_n)=\{\{x_{i},x_{i+1}\}\mid 1\leq i\leq n \text{ and } x_{n+1}=x_1\}$.

\begin{proposition}\label{path}
For $P_n$ and $C_n$, the following hold:
\begin{enumerate}
	\item $\displaystyle \reg(I(P_n))=\left\lfloor\frac{n+2}{3}\right\rfloor + 1$ for $n \geq 1$.
	
	\item $\reg(I(P_n))=\reg(I(P_{n-3})) + 1$ for $n \geq 4.$

	\item If $n \equiv 0,1 \pmod{3}$, then 
	$\reg(I(C_n))=\reg(I(C_n \setminus \{x\}))=\reg(I(C_n \setminus N[x])) + 1$
	for $n \neq 3,4$, and $\reg(I(C_n))=\reg(I(C_n\setminus \{x\}))=2$ for $n=3,4$.
	
	\item If $n \equiv 2 \pmod{3}$, then 
	$\reg(I(C_n))=\reg(I(C_n \setminus \{x\})) + 1 = \reg(I(C_n \setminus N[x])) + 1.$
\end{enumerate}
\end{proposition}
\begin{proof}
 The proof follows from \cite[Theorem 4.7, Theorem 5.2]{beyarslan1}.
\end{proof}
The next lemma, which is the well-known Depth Lemma, will play a crucial role in this paper.       
\begin{lemma} \label{depth-lemma}
	If $0 \rightarrow N \rightarrow M \rightarrow L \rightarrow 0$ is a short exact sequence of finitely generated modules over a local ring $R$, then
	\begin{enumerate}
		\item If $\operatorname{depth}(M) < \operatorname{depth}(L)$, then $\operatorname{depth}(N) = \operatorname{depth}(M)$.
		\item If $\operatorname{depth}(M) = \operatorname{depth}(L)$, then $\operatorname{depth}(N) \geq \operatorname{depth}(M)$.
		\item If $\operatorname{depth}(M) > \operatorname{depth}(L)$, then $\operatorname{depth}(N) = \operatorname{depth}(L) + 1$.
	\end{enumerate}
\end{lemma}

Polarization is a dominant technique to obtain a square-free monomial ideal from a given monomial ideal.
\begin{definition}
	Suppose that  $u=x^{a_1}_1x_2^{a_2} \cdots x^{a_n}_n $ is a monomial in $R$. Then we define the polarization of $u$ to be the square-free monomial
	\begin{equation*}
		\mathcal{P} (u)  =  x_{11}x_{12} \cdots   x_{1a_1}x_{21}x_{22} \cdots x_{2a_2} \cdots  x_{n1}x_{n2} \cdots x_{na_n}
	\end{equation*}
	in the polynomial ring $R^\mathcal{P}=\mathbb{K}[x_{ij}  \mid  1  \leq  i  \leq  n, 1  \leq  j  \leq  a_i ]$. If  $ I \subset  R $ is a monomial
	ideal with $\mathcal{G}(I)$ = $ \{ u_1,\ldots, u_m   \}$, the polarization of $ I $, denoted by $I^\mathcal{P}$ is defined as:
	\begin{equation*}
		I^\mathcal{P} = (\mathcal{P} (u_1), \ldots ,\mathcal{P}(u_m)).\end{equation*}
In particular, $I^\mathcal{P}$ is a square-free monomial ideal in the polynomial ring $R^\mathcal{P}$.
\end{definition}

\section{$k$-edge weighted VEW Graphs and vertex covers} \label{sec3}
In this section, we introduce the concept of strong vertex cover for a $k$-edge weighted VEW graph $D$. 

\begin{definition}
Let $D$ be a $k$-edge weighted VEW graph.	A vertex cover $ C $ of $ D $ is a subset of $  V(D)  $ such that if $ (x, y) \in E(D) ,$ then
	$ x \in C $ or  $ y \in C . $ A vertex cover $ C $ of $ D $ is minimal if each proper subset of $ C $ is not a
	vertex cover of $ D. $ We set $(C)$ to be the ideal generated by the variables in $C.$
\end{definition}
\begin{remark}\label{s.v.00}
Let $D$ be a $k$-edge weighted VEW graph with the underlying simple graph $G$. Then $ C $ is a minimal vertex cover of $ D $ if and only if $ C $ is a minimal vertex cover of $ G. $
\end{remark}
\begin{definition}
Let $D$ be a $k$-edge weighted VEW graph.	Let $ x $ be a vertex of  $ D, $ then the sets $ N_D^+ (x) =	\{y : (x, y) \in E(D)\}  $ and  $ N_D^- (x) = \{y : (y, x) \in E(D)\}  $ are called the out-neighbourhood and the in-neighbourhood of $ x $, respectively. Moreover, the neighbourhood of $ x $ is the set $ N_D(x) = N_D^+ (x)\cup N_D^- (x)$ and $N_D[x]=N_D(x) \cup \{x\}$.  Define $\deg_D(x) = |N_D(x)|$ for  $ x \in V(D) $.  A vertex $ x \in V(D) $ is called a source vertex if $N_D (x)= N_D^+ (x) .$ 
	A vertex $ x \in V(D) $ is called a sink vertex if $N_D (x)= N_D^- (x) .$  
\end{definition}

\begin{definition}\label{L1L2L3}
	Let $ C $ be a vertex cover of a $k$-edge weighted VEW graph $ D $ and  $ C^c $ is the complement of $  C ,$ i.e., $ C^c = V(D) \setminus C. $ Then    	we define \vspace*{0cm}\\
	\hspace*{3cm}$ L_1(C) = \{x \in C ~|~ N_D^+ (x) \cap C^c \neq        \emptyset \}, $ \vspace*{0.2cm}\\
	\hspace*{2.85cm} $L_2(C) = \{x \in C ~|~x\notin L_1(C) ~\mbox{and}~  N_D^- (x) \cap C^c \neq   \emptyset \}$ and  \vspace*{0.2cm}\\ 
	\hspace*{2.8cm}  $ L_3(C) = C \setminus (L_1(C) \cup L_2(C))$.
\end{definition}

 \begin{proposition}\label{L3}  
	Let $ C $ be a vertex cover of a $k$-edge weighted VEW graph $ D  $.  Then $ L_3(C) =\{x\in  C~|~N_D(x) \subset C\}$. Moreover, $ L_3(C) =\emptyset $ if  and  only  if $ C $ is  a minimal vertex cover of $ D.$
\end{proposition}

\begin{proof}
 If $ N_D(x) \nsubset C   $, then by Definition \ref{L1L2L3}, $x \in  L_1(C) \cup L_2(C) $.  Therefore, $$ L_3(C) =\{x\in  C~|~N_D(x) \subset C\}.$$ 
 \par Moreover, $y \in L_3(C)$, i.e., $N_D[y] \subset C$ if and only if $C\setminus \{y\}$ is a vertex cover of $D$. Hence the proof follows.
\end{proof}

	

	

\begin{remark}\label{rem-sink-source}
In a $k$-edge weighted VEW graph $D$, if $x$ is a source vertex, then $w(x) \leq k$, whereas if $x$ is a sink vertex, then $w(x) \geq k$.
\end{remark}

\begin{notation}
	Set $V^{*} = V^{*}(D) = \{v \in V(D) \mid N_D^{-}(v) \neq \emptyset  ~\mbox{and}~ w(v) \neq k\}$, and for $k=1$,     $V^*=V^{+} = V^{+}(D) = \{v \in V(D)|~ w(v) \neq 1\}$.
\end{notation}

\begin{definition}
 	A vertex cover $ C $ of $ D $ is strong if for each $ x \in L_3(C)  $, there is $(y, x) \in E(D)$
 	such that $ y   \in L_2(C) \cup L_3(C)$ with $ y \in   V^{*}(D)$.
 \end{definition}

 \begin{remark}\label{rem-3.9}
 	A vertex cover $ C $ of $ D $ is strong if for each $ x \in L_3(C)  $ there exist $ y \in N_D^-(x) \cap  V^{*}(D) \cap [L_2(C) \cup L_3(C)]$.       
 	A vertex cover $ C $ of $ D $ is not strong if and only if    there is $ x \in L_3(C)  $ such that $ y \in N_D^-(x) $ satisfies $ y   \in L_1(C) $, or   $ w(y) = k $, or $  N_D^-(y) = \phi $. 	
 \end{remark}

\begin{example}
	Consider the  $6$-edge weighted VEW bipartite graph $D_1$ and the $7$-edge weighted VEW odd cycle $D_2$  as in Figure \ref{fig.0}.
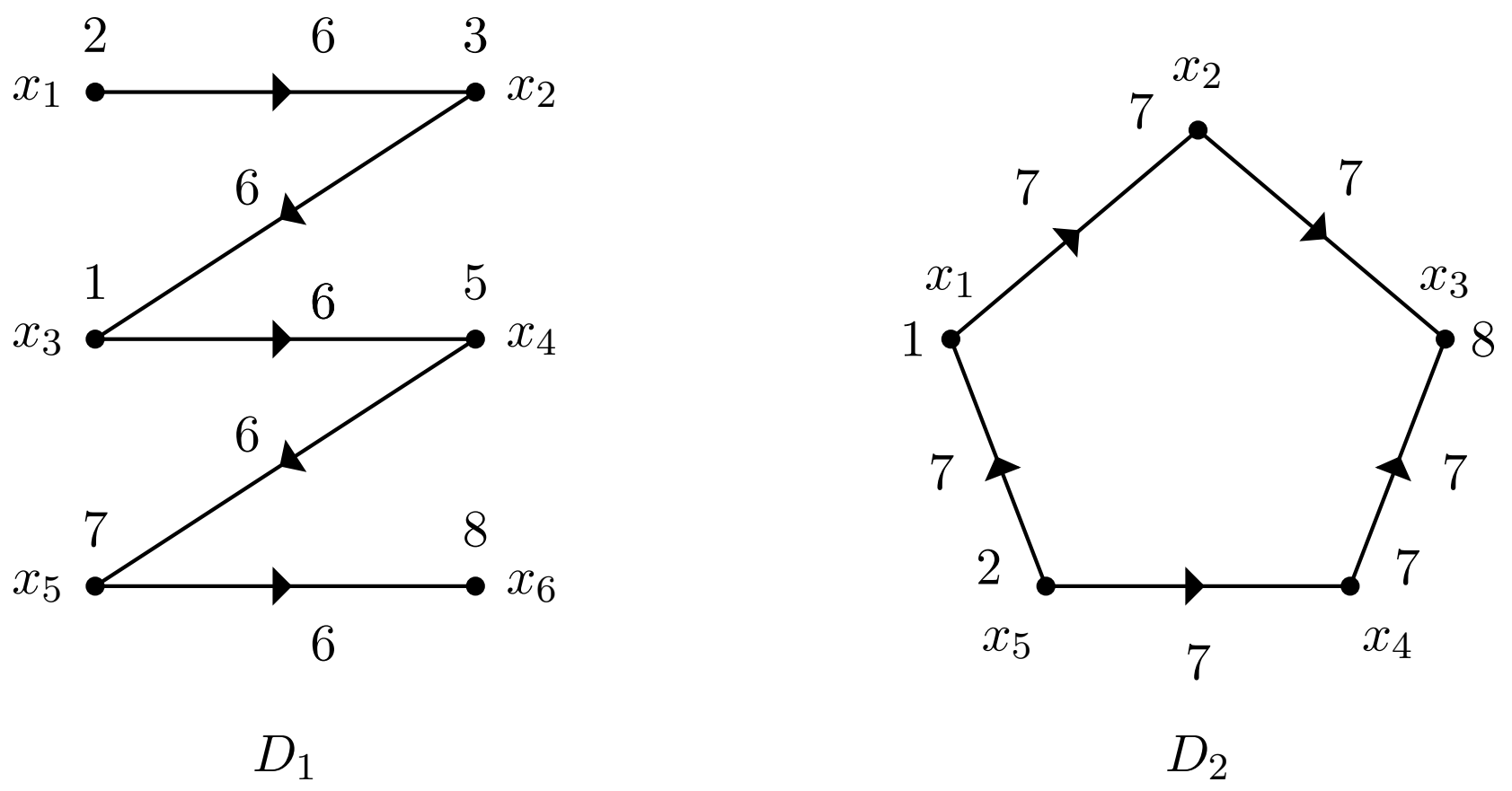
\begin{figure}[!ht]  
	\begin{tikzpicture}[scale=1.4]
	\begin{scope}[thick, every node/.style={sloped,allow upside down}] 
	\definecolor{ultramarine}{rgb}{0.07, 0.04, 0.56} 
	\definecolor{zaffre}{rgb}{0.0, 0.08, 0.66}

\node at (2.8,-0.3) {$x_5$};
\node at (4.8,-0.3) {$x_4$};
\node at (2.5,1.6) {$x_{1}$};
\node at (5.1,1.6) {$x_3$};
\node at (3.8,2.7) {$x_2$};

\node at (2.7,0.1) {$2$};  
\node at (4.9,0.1) {$7$};
\node at (2.3,1.3) {$1$};
\node at (5.3,1.3) {$8$};
\node at (3.5,2.5) {$7$}; 

\node at (2.9,2.1) {$7$};  
\node at (4.6,2.15) {$7$};
\node at (5.15,0.6) {$7$};
\node at (3.8,-0.4) {$7$};
\node at (2.45,0.6) {$7$};

\draw[fill] (3,0) --node {\midarrow}(4.6,0);
\draw[fill] (4.6,0) --node {\midarrow}(5.1,1.3);
\draw[fill] (3,0) --node {\midarrow}(2.5,1.3);
\draw[fill] (3.8,2.4) --node {\midarrow}(5.1,1.3);
\draw[fill] (2.5,1.3) --node {\midarrow}(3.8,2.4);

  \draw [fill] [fill] (3,0) circle [radius=0.04];
  \draw [fill] [fill] (4.6,0) circle [radius=0.04];
  \draw [fill] [fill] (5.1,1.3) circle [radius=0.04];
  \draw [fill] [fill] (3.8,2.4) circle [radius=0.04];
  \draw [fill] [fill] (2.5,1.3) circle [radius=0.04];

\node at (3.8,-0.9) {$D_2$};

\node at (0.3,0) {$x_6$};
\node at (0.3,1.3) {$x_4$};
\node at (0.3,2.6) {$x_{2}$};
\node at (-2.3,0) {$x_5$};
\node at (-2.3,1.3) {$x_3$};
\node at (-2.3,2.6) {$x_{1}$};

\node at (0,0.3) {$8$};
\node at (0,1.6) {$5$};
\node at (0,2.9) {$3$};
\node at (-2,0.3) {$7$};
\node at (-2,1.6) {$1$};
\node at (-2,2.9) {$2$};

\node at (-0.8,-0.3) {$6$};
\node at (-1.2,0.8) {$6$};
\node at (-0.8,1.5) {$6$};
\node at (-1.2,2.1) {$6$};
\node at (-0.8,2.9) {$6$};
\node at (-0.8,1.5) {$6$};

  \draw [fill] [fill] (-2,0) circle [radius=0.04];
  \draw [fill] [fill] (-0,0) circle [radius=0.04];
  \draw [fill] [fill] (-2,1.3) circle [radius=0.04];
  \draw [fill] [fill] (-0,1.3) circle [radius=0.04];
  \draw [fill] [fill] (-2,2.6) circle [radius=0.04];
  \draw [fill] [fill] (-0,2.6) circle [radius=0.04];

\draw[fill] (-2,2.6) --node {\midarrow}(-0,2.6);
\draw[fill] (-0,2.6) --node {\midarrow}(-2,1.3);
\draw[fill] (-2,1.3) --node {\midarrow}(-0,1.3);
\draw[fill] (-0,1.3) --node {\midarrow}(-2,0);
\draw[fill] (-2,0) --node {\midarrow}(-0,0);

\node at (-1,-0.9) {$D_1$};
	\end{scope}
	\end{tikzpicture}
		\caption{Two  $k$-edge weighted VEW graphs}\label{fig.0} 
\end{figure}
Note that $ V^*(D_1)= \{x_2,x_3,x_4,x_5,x_6 \}$ and $V^*(D_2)= \{x_1,x_3\}$. 
Consider one vertex cover $C_1=\{x_1,x_3,x_5,x_6\}$ of $D_1$. Here $L_1(C_1) = \{x_1,x_3\}$, $L_2(C_1) = \{x_5\}$ and $L_3(C_1) = \{x_6\}$. Since $x_5 \in N_{D_1}^-(x_6) \cap V^*({D_1})\cap L_2(C_1)$, $C_1$ is strong.
Now, we consider one vertex cover $C_2=\{x_1,x_2,x_3,x_4\}$ of $D_2$.  Here $L_1(C_2) = \emptyset$, $L_2(C_2) = \{x_1,x_4\}$ and $L_3(C_2) = \{x_2,x_3\}$. For $x_3 \in L_3(C_2)$, we have $ N_{D_2}^-(x_3) =\{x_2, x_4\}$ and $ \{x_2, x_4\} \cap V^*(D_2) = \phi $. Thus, $C_2$ is not a strong vertex cover of $D_2$.


\end{example}

\pagebreak
      
 \begin{remark}\label{VD}
If $C=V(D)$, then  $ L_1(C) = \emptyset,$ $ L_2(C) = \emptyset, $ and $ L_3(C) =  V(D)$.	Therefore,  $ V(D) $  is  strong if and only if for each $ x \in V(D)  $ there exist $ y \in N_D^-(x) \cap  V^{*}(D)$.  
\end{remark}

\begin{corollary}
If $ C $ is a minimal vertex cover of  $ D $, then $ C $ is strong.    
\end{corollary}

\begin{proof}
If $ C $ is minimal, then by Proposition \ref{L3}, $L_3(C) = \emptyset$. Hence it is strong.        	
\end{proof}

\begin{proposition}
	Let $ D $ be a $k$-edge weighted VEW graph. 
	Then $ V (D) $ is a strong vertex cover if
	and only if $ N_D^+(V^{*}(D)) = V(D) $. 	
\end{proposition}
\begin{proof}
	Let $ C = V (D). $ By Remark \ref{VD}, $ C $ is 
	strong  if and only if for each $ x \in L_3(C) = C $, there is $ (y,x) \in E(D) $ with $ y \in  V^{*}(D)$. Therefore, $ V (D) $ is a strong vertex
	cover if and only if $ N_D^+(V^{*}(D)) = V(D) $. 
\end{proof}

The following remark is a direct consequence of the above lemma.
\begin{remark}
If $ D $ is a $k$-edge weighted VEW graph  and contains a source vertex, then $ V (D) $    is not a strong vertex cover of $ D. $	
\end{remark}

\begin{proposition}
	Let $ D $ be a $k$-edge weighted VEW graph. Assume that
	$ w(v) \neq k $  for all $ v \in V(D) $.	Then $ V (D) $ is a strong vertex cover if and only if $D$ does not contain any source vertex. 	
\end{proposition}
 \begin{proof}
	By Remark \ref{VD}, $ V (D) $ is a strong vertex cover of $D$ if and only if for all $ x \in V (D) $ there is $ y \in V^*(D) $ such that $ (y,x) \in E(D). $ The later condition is equivalent to the condition ``$ D $ does not contain any source vertex".
 \end{proof}

\section{Edge ideals and their Primary decomposition}\label{sec4}

If $I$ is a monomial ideal of a polynomial ring $R$, there exists a  unique decomposition, $I=q_1 \cap \cdots \cap q_r$, where $q_1,\ldots,q_r$ are irreducible monomial ideals such that  $I \neq \bigcap_{i \neq j}q_i$  for each  $j=  1,\ldots,r$.   This is called the irredundant irreducible decomposition of $I$.  Moreover, $q_i$ is an irreducible monomial ideal if and only if $q_i= (x^{a_1}_{i_1},\ldots,x^{a_s}_{i_s})$ for some variables $x_{i_j}$.    Every irreducible decomposition is a primary decomposition as irreducible ideals are primary.  In this section, we define the edge ideal $I(D)$ of a $k$-edge weighted VEW graph $D$ and we characterize its irredundant irreducible decomposition.  Particularly, we prove that this decomposition is an irreducible primary decomposition, i.e, the radicals of the elements of the irredundant irreducible decomposition of $I(D)$ are different.

\begin{definition}
	Let $D$ be a $k$-edge weighted VEW graph.	Then the edge ideal of $D$ is defined as
	$$I(D)  =  \left(x_i^{\min\{k,w(x_i)\}} x_j^{\max\{k,w(x_j)\}}~|~(x_i,x_j) \in E(D) \right)      .$$	
	
\end{definition}

\begin{example}
Consider the  $6$-edge weighted VEW bipartite graph $D_1$ and the $7$-edge weighted VEW odd cycle $D_2$  as in Figure \ref{fig.0}.  Then $I(D_1) =(x_1^2x_2^6,x_2^3x_{3}^6,x_3x_4^6,x_4^5x_5^7,x_5^6x_6^8)$ 
 and  $I(D_2) =(x_1x_2^7,x_2^2x_{3}^7,x_4^7x_3^7,x_5^2x_4^7,x_5^2x_1^7)$.
 \end{example}

\begin{definition}\label{def:irreducible}
	Let $ C $ be a vertex cover of a $k$-edge weighted VEW graph $ D $. The irreducible ideal associated to $ C $ is the
	ideal
$$I_C  =  \left( \{x_i^{\min\{k, w(x_i)\}}~|~x_i \in L_1(C)   \}    \cup    \{x_j^{\max\{k, w(x_j)\}}~|~x_j \in  L_2(C) \cup L_3(C) \} \right).$$
\end{definition}

\begin{lemma}\label{edge}
	Let $ D $ be a $k$-edge weighted VEW graph. Then $ I(D) \subseteq  I_C $, for each vertex cover $ C $ of $ D $.
\end{lemma}

\begin{proof}
	Let $ I = I(D) $ and $ m \in \mathcal{G}(I) $. Then $ m = x_i^{\min\{k,w(x_i)\}} x_j^{\max\{k,w(x_j)\}}$ for some $(x_i,x_j) \in E(D) $. Since
	$ C $ is a vertex cover, $ x_i \in C $ or $ x_j \in C $.
	If $ x_j \in C $, then by Definition \ref{def:irreducible}, we have $ x_j^{k} \in I_C $ or $ x_j^{w(x_j)} \in I_C $. Thus
	$  x_j^{\max\{k,w(x_j)\}} \in I_C$ and so
	$ m = x_i^{\min\{k,w(x_i)\}} x_j^{\max\{k,w(x_j)\}} \in I_C$.
	Now, we assume that $x_j \notin C$. Then by Definition \ref{L1L2L3}, we have $x_i \in L_1(C) $ and
	$  x_i^{\min\{k,w(x_i)\}} \in I_C$.
Therefore
	$ m = x_i^{\min\{k,w(x_i)\}} x_j^{\max\{k,w(x_j)\}} \in I_C$.  Hence $I \subseteq I_C$.    
\end{proof}

\begin{lemma}\label{v.c}
	Let $ D $ be  a $k$-edge weighted VEW  graph  such  that $  I(D) \subseteq (x_{i_1}^{a_1},\ldots,x_{i_s}^{a_s}).$ Then $ \{{x_{i_1}},\ldots,{x_{i_s}}\} $ is a vertex cover of $ D $.   	
\end{lemma}

\begin{proof}
We take $ J = (x_{i_1}^{a_1},\ldots,x_{i_s}^{a_s}) $. If $ (a, b) \in E(D) $, then $a^{\min\{k,w(a)\}} b^{\max\{k,w(b)\}} \in I(D) \subseteq  J   $. Thus
$x^{a_j}_
{i_j}  ~|~  a^{\min\{k,w(a)\}} b^{\max\{k,w(b)\}}   $ for some $ 1 \leq j \leq s $. That means $x_
{i_j} \in \{a, b\}  $ and $ \{a, b\} \cap  \{  x_{i_1},\ldots,x_{i_s}\} \neq  \emptyset$.
Hence $ \{  x_{i_1},\ldots,x_{i_s}\} $ is a vertex cover of $ D $.  
\end{proof}

\begin{definition}
Let $ I $ be a monomial ideal. An irreducible monomial ideal q containing
$ I $ is called a minimal irreducible monomial ideal of $ I $ if for any irreducible monomial ideal
$ p $ such that $ I \subseteq p \subseteq q $ one has that $ p = q $.
\end{definition}

\begin{lemma}\label{minmax}
	Let $ D $ be a $k$-edge weighted VEW graph. Let $ J $ be a minimal irreducible monomial ideal of $I(D)$, where $\mathcal{G}(J) = \{x^{a_1}_
	{i_1} ,\ldots, x^{a_s}_
	{i_s}\}$. Then for any fixed $a_j$,  
    $$ {\min\{k,  w(x_{i_j}) \}} \leq  a_j \leq       {\max\{k, w(x_{i_j})\}}.$$        	
\end{lemma}

\begin{proof}
Here $I(D)  \subseteq   (x^{a_1}_
{i_1} ,\ldots, x^{a_s}_
{i_s})$. Let $M=\mathcal{G}(J)$. Fix some $j \in [s]$. First, suppose $   a_j  < {\min\{k,  w(x_{i_j}) \}}  .$   Then we take $J_1 = (M^{\prime})$, where $M^{\prime} = M \setminus \{x^{a_j}_
{i_j}  \}   \cup  \{x_{i_j}^{\min\{k,w(x_{i_j})\}}   \}$.    Note that $J_1 \subsetneq  J$.  Let $(a,b) \in E(D)$ and so $m = a^{\min\{k,w(a)\}} b^{\max\{k,w(b)\}} \in I(D) \subseteq  J   $.    Thus $  x^{a_r}_
{i_r}  ~|~   a^{\min\{k,w(a)\}} b^{\max\{k,w(b)\}}    $ for some $1 \leq r \leq s$.   If $r \neq j$, then $x^{a_r}_
{i_r} \in M^{\prime} $ and so $    a^{\min\{k,w(a)\}} b^{\max\{k,w(b)\}}  
\in J_1  $.
Now we assume $r = j$.  So $  x^{a_j}_
{i_j}  ~|~   a^{\min\{k,w(a)\}} b^{\max\{k,w(b)\}}    $. That means $  x_
{i_j} \in \{a,b\}$.
If $  x_
{i_j} = a$, then  $m=    x_{i_j}^{\min\{k,w(x_{i_j})\}} b^{\max\{k,w(b)\}}  \in J_1$ because $    x_{i_j}^{\min\{k,w(x_{i_j})\}}   \in J_1$.    If $  x_
{i_j} = b$, then  $m =     a^{\min\{k,w(a)\}} x_{i_j}^{\max\{k,w(x_{i_j})\}}  $.    Note that  $    x_{i_j}^{\max\{k,w(x_{i_j})\}}   \in J_1$, since ${\min\{k,w(x_{i_j})\}} \leq     {\max\{k,w(x_{i_j})\}} $. So  $    m   \in J_1$.  Therefore $ I(D) \subseteq J_1  \subsetneq  J $. Its a contradiction, since $J$ is minimal.  Hence  $    {\min\{k,  w(x_{i_j}) \}}  \leq  a_j.$ \par
Now, suppose $   a_j  > {\max\{k,  w(x_{i_j}) \}}  .$   Then we take $J_2 = (M^{\prime})$ where $M^{\prime} = M \setminus \{x^{a_j}_
{i_j}  \}$.    Note that $J_2 \subsetneq  J$.  Let $(a,b) \in E(D)$ and so $m = a^{\min\{k,w(a)\}} b^{\max\{k,w(b)\}} \in I(D)    $.  If  $ x_{i_j} \notin \{a,b\}$, then $  x^{a_l}_
{i_l}  ~|~   a^{\min\{k,w(a)\}} b^{\max\{k,w(b)\}}    $ for some $1 \leq l \leq s$ and $l \neq j$. So $m \in J_2$.    Note  that    $   a_j  > {\max\{k,  w(x_{i_j}) \}}  \geq {\min\{k,  w(x_{i_j}) \}}$.
If  $ x_{i_j} \in \{a,b\}$, then $  x^{a_j}_
{i_j}  ~\ndivides~   a^{\min\{k,w(a)\}} b^{\max\{k,w(b)\}}    $.
Thus $  x^{a_r}_
{i_r}  ~|~   a^{\min\{k,w(a)\}} b^{\max\{k,w(b)\}}    $ for some $1 \leq r \leq s$ and $r \neq j$. So $m \in J_2$. Therefore, $ I(D) \subseteq J_2  \subsetneq  J $. Its a contradiction, since $J$ is minimal. Hence $   a_j  \leq {\max\{k,  w(x_{i_j}) \}}  .$    
\end{proof}

\begin{lemma}\label{weights}
	Let $ D $ be a $k$-edge weighted VEW graph. Let $ J $ be a minimal irreducible monomial ideal of I(D), where $\mathcal{G}(J) = \{x^{a_1}_
	{i_1} ,\ldots, x^{a_s}_
	{i_s}\}$. Set $C =  \{x_{i_1} ,\ldots, x_{i_s}\}
	$.  Then 
	
	(1)  $a_j = {\min\{k,  w(x_{i_j}) \}}$   if $  x_{i_j} \in  L_1(C) $, 
	
	(2) $a_j =      {\max\{k, w(x_{i_j})\}}$ if $x_{i_j} \in  L_2(C) \cup L_3(C)  $.	
	
\end{lemma}

\begin{proof}
Here $I(D)  \subseteq   (x^{a_1}_
{i_1} ,\ldots, x^{a_s}_
{i_s})$. Since $\mathcal{G}(J)$ is minimal, we have $ x_{i_1} , \ldots , x_{i_s} $ are different. Note that by Lemma \ref{v.c}, $C$ is a vertex cover of $D$. 
	
	(1)  Let $  x_{i_j} \in  L_1(C) $. That means there is $x \in N_D^+(x_{i_j} )  \cap  C^c$.  Note that $(x_{i_j}, x  )  \in E(D)$. Here     $x_{i_j}^{\min\{k,w(x_{i_j})\}} x^{\max\{k,w(x)\}} \in I(D) \subseteq  J   $.    Since  $x \notin C$ and    $ x^{a_j}_
	{i_j}  ~|~       x_{i_j}^{\min\{k,w(x_{i_j})\}} x^{\max\{k,w(x)\}}    $, we have   $ x^{a_j}_
	{i_j}  ~|~       x_{i_j}^{\min\{k,w(x_{i_j})\}}    $.  So  $ {a_j}
 \leq       {\min\{k,w(x_{i_j})\}}    $.     
By Lemma \ref{minmax},  $ {a_j}
\geq       {\min\{k,w(x_{i_j})\}}    $. Therefore, $ {a_j}
=    {\min\{k,w(x_{i_j})\}}    $.

(2) Note that  $$I_C  =  ( \{x_{i_j}^{\min\{k, w(x_{i_j})\}}~|~x_{i_j} \in L_1(C)   \}    \cup    \{x_{i_j}^{\max\{k, w(x_{i_j})\}}~|~x_{i_j} \in  L_2(C) \cup L_3(C) \} ).$$

If $x_{i_j} \in L_1(C)  $, then by (1), $ {a_j}
=    {\min\{k,w(x_{i_j})\}}    $.
If  $x_{i_j} \in  L_2(C) \cup L_3(C)  $,  then by Lemma \ref{minmax},  $ {a_j}
\leq       \max\{k,w(x_{i_j})\}   $. Therefore,    $ I_C \subseteq J$.  If  $ {a_j}
<       {\max\{k,w(x_{i_j})\}}    $ for some $j \in [s]$, then $ I_C \subsetneq J$.
Thus, by Lemma \ref{edge}, $I(D)  \subseteq I_C \subsetneq J$. It's a contradiction because $J$ is minimal. Hence, $a_j ={\max\{k, w(x_{i_j})\}}$ if $x_{i_j} \in  L_2(C) \cup L_3(C)  $.    
\end{proof}

\begin{theorem}\label{minimal-j}
The following conditions are equivalent:

$ (1) $ $ J $ is a minimal irreducible monomial ideal of $ I(D) $.

$ (2) $ There is a strong vertex cover $ C $ of $ D $ such that $ J = I_C $.
\end{theorem}

\begin{proof}
$ 1) \Rightarrow  2)$  Let $ J $ be a minimal irreducible monomial ideal of $ I(D) $. By  Lemma \ref{weights}, $J =I_C= (x^{a_1}_
{i_1} ,\ldots, x^{a_s}_
{i_s})$ for some vertex cover $C = \{x_
{i_1} ,\ldots, x_
{i_s}\}$ of $D$. Suppose that $C$ is not strong. Then by Remark \ref{rem-3.9}, there is $ x \in L_3(C) $ such that if $ (y, x) \in E(D) $, then
$ y   \in L_1(C) $, or   $ w(y) = k $, or $  N_D^-(y) = \phi $. 	 Without loss of generality, let $ x = x_{i_1} $  and set $ J^{\prime} $ as the monomial ideal
with $ \mathcal{G}(J^{\prime}) = ( x_{i_2}^{a_2}, \ldots, x_{i_s}^{a_s} )$,  where $a_j = {\min\{k,  w(x_{i_j}) \}}$   if $  x_{i_j} \in  L_1(C) $ and $a_j =      {\max\{k, w(x_{i_j})\}}$ if $x_{i_j} \in  L_2(C) \cup L_3(C)  $ for $ 2 \leq j \leq s$.

  Let $ (z_1, z_2) \in E(D) $. If $ x_{i_j}^{a_j}
  ~|~  {z_1}^{\min\{k,w(z_1)\}}{z_2}^{\max\{k,w(z_2)\}} $
 for some $ 2 \leq j \leq s $,
then  $   {z_1}^{\min\{k,w(z_1)\}}{z_2}^{\max\{k,w(z_2)\}} \in  J^{\prime}$.
 Now, we assume $ x_{i_j}^{a_j}
 \ndivides  {z_1}^{\min\{k,w(z_1)\}}{z_2}^{\max\{k,w(z_2)\}} $ for each $ 2 \leq j \leq s $.  Note that $a_j \in \{\min\{k,w(x_{i_j})\}, \max\{k,w(x_{i_j})\}    \} $.
If  $ z_2 \in  \{x_{i_2},  \ldots    ,x_{i_s}\}  $, then  $ x_{i_j}^{a_j}
\divides  {z_1}^{\min\{k,w(z_1)\}}{z_2}^{\max\{k,w(z_2)\}} $ for some $ 2 \leq j \leq s $, which contradicts our assumption.  Therefore
 $ z_2 \notin  \{x_{i_2},  \ldots    ,x_{i_s}\}  $ and so  $ z_2 = x_{i_1} $.  Since $ (z_1, x_{i_1}) \in E(D) $ and $C$ is not strong, by Remark \ref{rem-3.9}, $ w(z_1) = k $, or $ z_1 \in L_1(C)  $, or $  N_D^-(z_1) = \phi $. In all cases $ {z_1}^{\min\{k,w(z_1)\}} \in I_C = J $.  Since
 $ z_1 \neq x_{i_1} $, we get $ {z_1}^{\min\{k,w(z_1)\}} \in J^{\prime} $. This implies $ {z_1}^{\min\{k,w(z_1)\}}{z_2}^{\max\{k,w(z_2)\}}
 \in  J^{\prime} $. Hence, $ I(D) \subseteq J^{\prime} $. This is a contradiction, since $ J $ is minimal
 and $ J^{\prime}   \subsetneq J $. Therefore, C is strong vertex cover of $D$.       

$ 2) \Rightarrow  1)$  Assume that there is a strong vertex cover $ C $ of $ D $ such that $ J = I_C $.
By definition $ J = I_C $ is a monomial irreducible ideal. By Lemma \ref{edge},
$ I(D)  \subseteq  J $. Now, suppose $ I(D) \subseteq J^{\prime} \subseteq J $, where $ J^{\prime} $ is a monomial irreducible ideal. We
can assume $ \mathcal{G}(J^{\prime}) = ( x_{j_1}^{b_1}, \ldots, x_{j_s}^{b_s} )$.

If $ x \in L_1(C) $, then there is $ (x, y) \in E(D) $ with $ y \notin C $. Then $ {x}^{\min\{k,w(x)\}}{y}^{\max\{k,w(y)\}} \in I(D) $ and $ y^r  \notin J $
for each $ r \in \mathbb{N} $. Thus $ y^r  \notin J^{\prime}$ for each $ r $ implying $y \notin \{ x_{j_1}, \ldots, x_{j_s} \} $ . 
 Since $ {x}^{\min\{k,w(x)\}}{y}^{\max\{k,w(y)\}} \in   J^{\prime},$ 
 $ x_{j_i}^{b_i}
 ~|~  {x}^{\min\{k,w(x)\}}{y}^{\max\{k,w(y)\}} $ for some $ 1 \leq i \leq s $.  Note that $ 
x_{j_i}^{b_i}  ~|~ {x}^{\min\{k,w(x)\}}$ and so  $ {x}^{\min\{k,w(x)\}}   \in  J^{\prime}$.

If $ x \in L_2(C) $, then there is $ (y, x) \in E(D) $ with $ y \notin C $.  Thus $y^{\min\{k,w(y)\}} \notin ( x_{j_1}^{b_1}, \ldots, x_{j_s}^{b_s} ) $, since $y \notin \{ x_{j_1}, \ldots, x_{j_s} \} $. Also, $ {y}^{\min\{k,w(y)\}}{x}^{\max\{k,w(x)\}} \in I(D) \subseteq  J^{\prime}$. That means $ x_{j_i}^{b_i}
~|~\\  {y}^{\min\{k,w(y)\}}{x}^{\max\{k,w(x)\}} $ for some $ 1 \leq i \leq s $.   Therefore,  $ x_{j_i}^{b_i}
~|~  {x}^{\max\{k,w(x)\}} $  implies $ {x}^{\max\{k,w(x)\}}\\ \in J^{\prime}. $

Finally, we assume that  $ x  \in L_3(C) $. Since $ C $ is a strong,  there is $ (y, x) \in E(D) $, where $ y \in L_2(C) \cup L_3(C) $ and $ w(y) \neq k $.  
 So, $ {y}^{\min\{k,w(y)\}}{x}^{\max\{k,w(x)\}} \in I(D) \subseteq  J^{\prime} $ implies $ x_{j_i}^{b_i}
~|~  {y}^{\min\{k,w(y)\}}{x}^{\max\{k,w(x)\}} $ for some $ 1 \leq i \leq s $. Since $ w(y) \neq k $, we have $ {\min\{k,w(y)\}} < {\max\{k,w(y)\}} $. That means $ {y}^{\min\{k,w(y)\}} \notin J = I_C $, because $ y  \in L_2(C) \cup L_3(C) $. This implies $ {y}^{\min\{k,w(y)\}} \notin J^{\prime} $, and   so,  $ x_{j_i}^{b_i}
~|~  {x}^{\max\{k,w(x)\}} $, then $ {x}^{\max\{k,w(x)\}} \in J^{\prime}. $
Hence, $  J = I_C \subseteq  J^{\prime}$. Therefore, $ J $ is a minimal monomial irreducible ideal of $ I(D) $.
\end{proof} 

\begin{theorem} \label{irreducible}
If $ C_s $ is the set of strong vertex covers of $ D $, then the irredundant irreducible
decomposition of $ I(D) $ is given by $ I(D) =  \bigcap_{C \in C_s}  I_C. $

\end{theorem}

\begin{proof}

By \cite[Theorem 1.3.1]{hh11}, there exists a unique irredundant irreducible decomposition $I(D)=\bigcap_{i=1}^{t} I_i$. Suppose there exists an irreducible ideal $I'_j$ such that
$I(D)\subseteq I'_j \subseteq I_j$ for some $j\in \{1,\ldots,t\}$. Then $I(D)=\left(\bigcap_{i\ne j} I_i\right)\cap I'_j$ is an irreducible decomposition. Also, this decomposition is irredundant, since
\(\bigcap_{i=1}^{t} I_i\) is irredundant. Thus $I'_j = I_j$, since the irredundant irreducible decomposition is unique. Hence, $I_1,\ldots,I_t$ are minimal irreducible ideals of $I(D)$. Then, by Theorem \ref{minimal-j}, $\{I_1,\ldots,I_t\} \subseteq \{I_C \mid C\in \mathcal{C}_s\}$. Now, suppose there exists $C\in \mathcal{C}_s$ such that
$I_C \notin \{I_1,\ldots,I_t\}$. Then there exists $x^{\alpha_{i}}_{j_i}\in I_i \setminus I_C$ for each $i\in \{1,\ldots,t\}$.
Consequently, $q=\operatorname{lcm}(x^{\alpha_1}_{j_1},\ldots,x^{\alpha_{t}}_{j_t})
\in \bigcap_{i=1}^{t} I_i = I(D) \subseteq I_C.$ Moreover, if $C=\{x_{i_1},\ldots,x_{i_l}\}$, then
$I_C=(x_{i_1}^{\gamma_1},\ldots,x_{i_l}^{\gamma_l})$,
where $\gamma_j \in \{\min\{k,w(x_{i_j})\}, \max\{k,w(x_{i,j})\}\}$. Hence, there exists $j\in \{1,\ldots,l\}$ such that $x_{i_j}^{\gamma_j} \mid q$.
Thus, there exists $1\le u\le t$ such that
$x_{i_j}^{\gamma_j} \mid x_{j_u}^{\alpha_u}$,
which is a contradiction, since $x_{j_u}^{\alpha_u} \notin I_C$. Therefore,
$\bigcap_{i=1}^{t} I_i = \bigcap_{C\in \mathcal{C}_s} I_C$
is the irredundant irreducible decomposition of $I(D)$.
\end{proof}

\begin{remark} \label{prime}
 Assuming $C_1,\ldots,C_r$ are the strong vertex covers of $D$,  by Theorem \ref{irreducible}, $I_{C_1} \cap \cdots \cap I_{C_r}$ is the irredundant irreducible decomposition of $I(D)$. If $P_i = \rad(I_{C_i} )$,
then $P_i = (C_i)$. So, $P_i \neq P_j$ for $1 \leq i < j \leq   r$. Therefore, $I_{C_1} \cap \cdots \cap I_{C_r}$ is an irredundant
primary decomposition of $I(D)$  and $\Ass(I(D)) = \{P_1, \ldots, P_r\}$. 
\end{remark}
\begin{example}
Consider the $6$-edge weighted VEW odd cycle $D$  as in Figure \ref{fig.2}.  Then, $I(D) =(x_1^6x_2^7,x_2^6x_{3}^6,x_4^6x_3^6,x_5^2x_4^9,x_5^2x_1^8)$.   
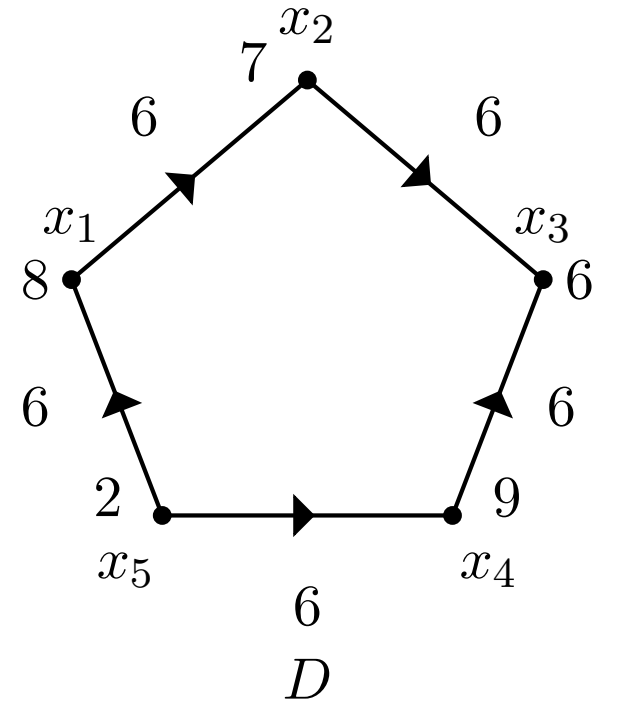
\begin{figure}[!ht]  
	\begin{tikzpicture}[scale=1.2]
	\begin{scope}[thick, every node/.style={sloped,allow upside down}] 
	\definecolor{ultramarine}{rgb}{0.07, 0.04, 0.56} 
	\definecolor{zaffre}{rgb}{0.0, 0.08, 0.66}

\node at (2.8,-0.3) {$x_5$};
\node at (4.8,-0.3) {$x_4$};
\node at (2.5,1.6) {$x_{1}$};
\node at (5.1,1.6) {$x_3$};
\node at (3.8,2.7) {$x_2$};

\node at (2.7,0.1) {$2$};  
\node at (4.9,0.1) {$9$};
\node at (2.3,1.3) {$8$};
\node at (5.3,1.3) {$6$};
\node at (3.5,2.5) {$7$}; 

\node at (2.9,2.2) {$6$};  
\node at (4.8,2.2) {$6$};
\node at (5.2,0.6) {$6$};
\node at (3.8,-0.5) {$6$};
\node at (2.3,0.6) {$6$};

\draw[fill] (3,0) --node {\midarrow}(4.6,0);
\draw[fill] (4.6,0) --node {\midarrow}(5.1,1.3);
\draw[fill] (3,0) --node {\midarrow}(2.5,1.3);
\draw[fill] (3.8,2.4) --node {\midarrow}(5.1,1.3);
\draw[fill] (2.5,1.3) --node {\midarrow}(3.8,2.4);

  \draw [fill] [fill] (3,0) circle [radius=0.04];
  \draw [fill] [fill] (4.6,0) circle [radius=0.04];
  \draw [fill] [fill] (5.1,1.3) circle [radius=0.04];
  \draw [fill] [fill] (3.8,2.4) circle [radius=0.04];
  \draw [fill] [fill] (2.5,1.3) circle [radius=0.04];

	\node at (3.8,-0.9) {$D$};

	\end{scope}
	\end{tikzpicture}
		\caption{A $6$-edge weighted VEW cycle of length 5}\label{fig.2} 
\end{figure}
 By Theorem \ref{irreducible},  the irredundant irreducible decomposition of $I(D)$ is given as
\begin{align*}
I(D) 
&=  (x_1^8,x_2^6,x_4^6) \cap 
 (x_1^6,x_3^6,x_4^9) \cap  (x_1^6,x_3^6,x_5^2) \cap  (x_2^7,x_3^6, x_5^2)\\
&\hspace{0.4cm}  \cap  (x_2^6,x_4^6,x_5^2)  \cap  (x_1^8, x_2^7, x_3^6, x_4^9),    
	\end{align*}
 and $\Ass(I(D)) = \{(x_1,x_2,x_4) , 
 (x_1,x_3,x_4),  (x_1,x_3,x_5),  (x_2,x_3, x_5),  (x_2,x_4,x_5) ,  (x_1, x_2, x_3,\\ x_4)    \}.$        
\end{example}
\section{Cohen-Macaulayness of  $k$-edge weighted VEW  graphs}\label{sec5}
In this section, we give criteria for the unmixedness and Cohen--Macaulayness of edge ideals of certain classes of $k$-edge weighted VEW graphs. Recall that an ideal $I$ is called unmixed if all its associated primes have the same height. It is well known that if $R/I$ is Cohen--Macaulay, then $I$ is unmixed, but the converse may not hold in general. The following theorem shows that the converse holds for forest.

\begin{theorem}\label{cm-forest}
Let $D$ be a $k$-edge weighted VEW forest without isolated vertices, and let $G$ be its underlying forest. Then the following statements are equivalent:
\begin{enumerate}
    \item $D$ is Cohen--Macaulay.
    
    \item $I(D)$ is unmixed.
    
    \item $G$ admits a perfect matching $\{\{x_1,y_1\}, \{x_2,y_2\}, \ldots, \{x_r,y_r\}\}$ such that 
    $\deg_D(y_i)=1$ for all $i=1,2,\ldots,r$, and for each $i$ with $(x_i,y_i) \in E(D)$, if $x_i$ is not a source vertex then $w(x_i)=k$.
\end{enumerate}
\end{theorem}
\begin{proof}
	$(1) \implies (2):$ Follows from \cite[Corollary 3.1.17]{v15}. \\
	$(2) \implies (3):$ By \cite[Theorem 2.4]{v90} and \cite[Corollary 2.5]{v90}, $|V(G)| = 2r$, and $G$ has a perfect matching $\{x_1, y_1\}, \dots, \{x_r, y_r\}$ with $\deg_D(y_i) = 1$ for $i = 1, \dots, r$.
Let $D'$ be an oriented graph with vertex set $V(D') = \{x_1, \ldots, x_r\}$ and edge set $E(D')=\{(x_i, x_j) \mid (x_i, x_j) \in E(D)\}$. Note that $D'$ is acyclic. Then, by \cite[Theorem 16.3(4),p.200]{h69}, there exists a linear ordering of the vertex set $V(D')$ such that all the edges of $D'$ are of the form $(x_i, x_j)$ with $i < j$. There is nothing to prove if $(y_i, x_i) \in E(D)$ for $i = 1, \ldots, r$. Now, assume $(x_l, y_l) \in E(D)$ for some $l$ and $x_l$ is a non-source vertex. Then we have to show that $w(x_l) = k$.  Suppose $w(x_l) \neq k$. Suppose $A = \{x_1, \ldots, x_r\}$ and set $C = (A\setminus N^{-}_{D'}(x_l)) \cup \{y_i \mid x_i \in N^{-}_{D'}(x_l)\} \cup \{y_l\},$ where $N^{-}_{D'}(x_l)=\{y \in V(D') \mid (y,x_l) \in E(D')\}$. Note that $N^{-}_{D'}(x_l)$ is an independent set of $G$. Thus, $C$ is a vertex cover of $G$ with $r+1$ elements. Now, we show that $C$ is a strong cover of $D$. Note that $x_l$ cannot be a source of $D'$. Thus, the set $N^{-}_{D'}(x_l) \neq \emptyset $. This implies that $x_l \not \in L_3(C)$. Therefore, $L_3(C) = \{y_l\}$ because $L_3(C) \subset \{x_l, y_l\} \subset C$. Also, $N_D^+ \cap C^c = \emptyset$ and $N_D^- \cap C^c \neq \emptyset$. Thus, $x_l \in L_2(C)$. Hence, $(x_l, y_l) \in E(D)$ with $x_l \in L_2(C)$ and $w(x_l) \neq k$. This means that $C$ is a strong cover of $D$. Thus, by Remark \ref{prime}, $P =(C)$ is an associated prime of $I(D)$ with $r+1$ elements, which is a contradiction because $I(D)$ is an unmixed ideal of height $r$. \\
	$(3) \implies (1):$
	We argue by induction on $r=\frac{\mid V(D) \mid}{2}$. The base case $r = 1$ is clear, as $I(D)$ is a principal ideal, thus Cohen–Macaulay. Let $D'$ be the induced oriented subgraph of $D$ on the vertex set $A = \{x_1, \ldots, x_r\}$.  By \cite[ Theorem 16.3(4), p. 200]{h69}, there exists a linear ordering of the vertex set $V(D')$ such that all the edges of $D'$ are of the form $(x_i, x_j)$ with $i < j$. Set $R = \mathbb{K}[x_1, \ldots, x_r, y_1, \ldots, y_r]$. Now, we consider two possible cases. \\
	\noindent
	\textbf{Case (1)} Suppose $(x_r, y_r) \in E(D)$. Then $w(x_r) = k$ because $x_r$ is a non-source vertex in $D$. 
   
    Therefore, $x_r^{\min\{w(x_r),k\}}y_r^{\max\{w(y_r),k\}}=x_r^{k}y_r^{w(y_r)} \in I(D)$ because $y_r$ is a sink vertex.  Note that
	\begin{align*}
		I(D): x_r^{k} =& \left((z^{\min\{w(z),k\}} \mid z \in N_D(x_r)\setminus\{y_r\}), y_r^{w(y_r)}, I(D\setminus N_D(x_r))\right) \quad \text{and} \quad \\
		(I(D), x_r^{k}) = &(x_r^{k}, I(D\setminus\{x_r\})).
	\end{align*}
Thus, by applying induction hypothesis to $I(D\setminus N_D(x_r))$ and $I(D\setminus\{x_r\})$, we have $I(D): x_r^{k}$ and $(I(D), x_r^{k})$ are Cohen–Macaulay of dimension $r$. Now, let us take an exact sequence
	\begin{align}\label{unmixed-forest2}
		0 \rightarrow \frac{R}{I(D): x_r^{k}}(-k) \xrightarrow{ .x_r^{k} } \frac{R}{I(D)} \rightarrow \frac{R}{(I(D),x_r^{k})} \rightarrow 0 . 
	\end{align}
	Note that $I(D)$ has height $r$. Thus, by applying Lemma \ref{depth-lemma} to exact sequence \eqref{unmixed-forest2}, we have $I(D)$ is Cohen–Macaulay.  \\
	\noindent
	\textbf{Case (2)} Suppose $(y_r, x_r) \in E(D)$. Then $x_r$ is a sink of $D$. Note that
	\begin{align*}
		I(D): x_r^{w(x_r)} =& \left((z^{\min\{w(z),k\}} \mid z \in N_D(x_r)), I(D\setminus N_D(x_r))\right) \quad \text{and} \\\quad (I(D), x_r^{w(x_r)}) =& \left(x_r^{w(x_r)}, I(D\setminus\{x_r\})\right).
	\end{align*}
Thus, by applying induction hypothesis to $I(D\setminus N_D(x_r))$ and $I(D\setminus \{x_r\})$, we have $I(D): x_r^{w(x_r)}$ and $(I(D), x_r^{w(x_r)})$ are Cohen–Macaulay of dimension $r$. Now, let us take an exact sequence
	\begin{align}\label{unmixed-forest}
		0 \rightarrow \frac{R}{I(D): x_r^{w(x_r)}}(-w(x_r)) \xrightarrow{ .x_r^{w(x_r)} }\frac{R}{I(D)} \rightarrow \frac{R}{(I(D), x_r^{w(x_r)})} \rightarrow 0 . 
	\end{align}
	Note that $I(D)$ has height $r$. Thus, by applying Lemma \ref{depth-lemma} to exact sequence \eqref{unmixed-forest}, we have $I(D)$ is Cohen–Macaulay. 
\end{proof}
While the equivalence in Theorem~\ref{cm-forest} holds when the underlying graph $G$ is a forest, the following theorem shows that a similar equivalence also holds for certain $k$-edge weighted VEW graphs whose underlying graphs admit a perfect matching.
\begin{theorem}\label{cm-graph}
	Let $D$ be a $k$-edge weighted VEW graph without isolated vertices and the underlying graph $G$ has a perfect matching $\{x_1,y_1\}, \{x_2,y_2\}, \ldots, \{x_r,y_r\}$ such that $\deg_D(y_i)=1$ for $i=1,2,\ldots,r$. 
	Then the following conditions are equivalent: 
	\begin{enumerate}
		\item $D$ is Cohen-Macaulay. 
		\item $I(D)$ is unmixed.
		\item For each $i$ with $(x_i,y_i) \in E(D)$, if $x_i$ is a non-source vertex, then $w(x_i)=k$.
	\end{enumerate} 
\end{theorem}
\begin{proof}
	$(1) \implies (2):$ Follows from \cite[Corollary 3.1.17]{v15}. \\
	$(2) \implies (3):$ Suppose $I(D)$ is unmixed. There is nothing to prove if $(y_i,x_i) \in E(D)$ for $i=1,\ldots,r$. Now, suppose that $(x_l, y_l) \in E(D)$ for some $l$ and $x_l$ is a non-source vertex. We need to show that $w(x_l)=k$. Suppose $w(x_l) \neq k$. Then, there exist at least one edge $(x_s,x_l) \in E(D)$ because $y_i$'s are leaf vertices in $G$ and $x_l$ is a non-source vertex. Note that $C = \left(\{x_1, \ldots, x_r\} \backslash \{x_s\}\right) \cup \{y_l, y_s\}$ is a vertex cover of $G$ which is not minimal. Also, $L_3(C) = \{y_l\}$ and $x_l \notin L_1(C)$. This implies that $x_l \in L_2(C)$. Thus, $C$ is a strong vertex cover of $D$ because $(x_l, y_l) \in E(D)$ and $w(x_l) \neq k$. Therefore, $I(D)$ is not unmixed, which is a contradiction. Hence, for each $i$ with $(x_i,y_i) \in E(D)$, if  $x_i$ is not a source vertex, then $w(x_i)=k$. 
	
	$(3)\implies (1)$. We will prove this implication using polarization. Note that $y_i's$ are free variables in $I(D)$ for all $1 \leq i \leq r$. Thus, by \cite[Lemma 2.8]{H.T-2018}, suppose the exponents of $y_i's$ is $1$ in $I(D)$ for all $i = 1, \ldots, r$. Let $D[A]$ be the induced oriented subgraph of $D$ on the vertex set $A = \{x_1, \ldots, x_r\}$. By re-indexing if necessary, we may assume that for some $l \leq r$, $w(x_i) > k$ for $1 \leq i \leq l$ and $w(x_i) \leq k$ for $l < i \leq r$. Thus, we have $(y_i, x_i) \in E(D)$ for $i = 1, \ldots, l$ because of Assumption (3) and by Remark \ref{rem-sink-source}. Now, define for all $l < i \leq r$,
	\begin{align*}
		a_i= & k \text{   if  } (y_i,x_i) \in E(D)  \\
		= &w(x_i) \text{  if  } (x_i,y_i) \in E(D).
	\end{align*}
	Then,
	$ I(D) = (y_1x^{w(x_1)}_1, \ldots, y_lx^{w(x_l)}_l, x_{l+1}^{a_{l+1}}y_{l+1}, \ldots, x_r^{a_r}y_r) + I(D[A]).$ Suppose $H$ is a copy of $D[A]$ on the vertex set $B = \{u_1, \ldots, u_r\}$ obtained by replacing $x_i$ with $u_i$ for $i = 1, \ldots, r,$ and let $S = \mathbb{K}[u_1, \ldots, u_r]$. Note that, in $H$, $w(u_i)=w(x_i)$ for all $i=1, \ldots, r$.  Consider the ideal 
	$$J = I(H) + (u_1^{w(u_1)+1}, \ldots, u_l^{w(u_l)+1},u_{l+1}^{a_{l+1}+1},\ldots,u_{r}^{a_r+1}) \subseteq S.$$ Note that in $J^{\PP}$, the variables $u_1, \ldots, u_r$ are replaced by $u_{1,1}, \ldots, u_{r,1}$. Also, for $i = 1, \ldots, l$, we have $(u_i^{w(u_i)+1})^{\PP} = u_{i,1} \ldots u_{i,w(u_i)+1}$ and for $i = l+1, \ldots, r$, we have $(u_i^{a_i+1})^{\PP} = u_{i,1} \ldots u_{i,a_i+1}$. On the other hand, in $I(D)^{\PP}$, the variables $x_1, \ldots, x_r$ are replaced by $x_{1,1}, \ldots, x_{r,1}$. Also, for $i = 1, \ldots, l$, we have
	$(y_ix_i^{w(x_i)})^{\PP} = y_{i,1}x_{i,1} \ldots x_{i,w(x_i)}$ and for $i = l+1, \ldots, r$, we have $(y_ix_i^{a_i})^{\PP} = y_{i,1}x_{i,1} \ldots x_{i,a_i}$. Now, consider $\psi : S^{\PP}\rightarrow R^{\PP}$ given by 
	\begin{align*}
		u_{i,j} \rightarrow & x_{i,j} \text{ for } 1 \leq i \leq l \text{ and } 1 \leq j \leq w(x_i), \\
		u_{i,j} \rightarrow & x_{i,j} \text{ for } l+1 \leq i \leq r \text{ and } 1 \leq j \leq a_i, \\
		u_{i, w(u_i)+1} \rightarrow & y_{i,1} \text{ for } 1 \leq i \leq l,~~ u_{i, a_i+1} \rightarrow  y_{i,1} \text{ for } l+1 \leq i \leq r.
	\end{align*}
	Then $I(D)^{\PP} = \psi(J^{\PP})$. Note that $S/J$ is an Artinian ring. This implies that $S/J$ is Cohen–Macaulay. Thus, by \cite[Proposition 2.8]{f}, $S^{\PP}/J^{\PP}$ is Cohen–Macaulay. Therefore, $R^{\PP}/I(D)^{\PP}$ is Cohen–Macaulay. Again by \cite[Proposition 2.8]{f}, $R/I(D)$ is Cohen-Macaulay.
\end{proof}

\section{Regularity of $k$-edge weighted VEW graphs} \label{sec6}
In this section, we compute the regularity of certain classes of  VEW graphs by connecting their polarized edge ideal with the labeled hypergraph and using the technique of Betti splitting. The concept of a labeled hypergraph is introduced in
\cite{McCullough}. Let $D$ be a VEW  graph.  The following result from \cite{McCullough} relates the regularity and projective dimension of $R/I(D)$
to the graph-theoretical invariants of the labeled hypergraph $H(D).$  
\begin{proposition}\cite[Proposition 4.1]{McCullough}\label{labeled}
Let D be a VEW  graph with associated labeled hypergraph $H = H(I(D)) = (V(H), X(H), \mathcal{E}(H)) $. If
	$\{i\} \in \mathcal{E}(H)$ for all $i \in [m] = V (H)$, then $\reg(R/I(D)) = |X(H)| - |V (H)|$ and $\pd(R/I(D)) = |V (H)|$.
\end{proposition}

Our next result is an application of the above proposition.

\begin{theorem}\label{atmost-one}
Let $ D $ be a VEW graph with the property that there is at most one edge oriented into each vertex and suppose that for all non-leaf, non-source vertices, $  x_j, $ either $ w(x_j) \neq k $ or the unique edge $ (x_i, x_j ) $ into the vertex $ x_j $ has the property that $ x_i $ is a leaf.
Then

 $\hspace*{2cm}\displaystyle{  \reg(R/I(D)) = \sum_{\substack{x_i \in V(D)\\~\mbox{and}~N^{-}_D(x_i) \neq \emptyset}} {{\max\{k,{w(x_{i})}\}}} + \sum_{\substack{x_i \in V(D)\\~\mbox{and}~N^{-}_D(x_i) = \emptyset}} {{{w_{i}}}}-|E(D)|}  ,$
\vspace*{0.2cm}\\
$\hspace*{2.1cm}\displaystyle{  \pd(R/I(D)) = |E(D)|}.$ 
\end{theorem}

\begin{proof}
Let $ w_i=w(x_i)  $ for every $x_i \in V(D).$	Let  $m_1, \ldots, m_r$ be the minimal
	generators of $I(D)^\mathcal{P}$, the polarization of $I(D)$. Let $ H(I(D)) = (V(H), X(H), \mathcal{E}(H)) $ be the labeled hypergraph associated to $I(D)$. For all $1 \leq  s \leq r$,  let  $\displaystyle{ m_s= \prod_{l=1}^{min\{k,{w_{i}}\}}x_{{i}l}\prod_{t=1}^{max\{k,{w_{j}}\}}x_{{j}t}}, $  i.e., a minimal generator of $ I(D)^\mathcal{P} $ corresponding to some edge $(x_{i},x_{j})$ of $D$.  First, we assume  $ x_{j} $ is a leaf in  $ D $. Then $ m_s $
	is the only minimal generator of $ I(D)^{\mathcal{P}} $
	which is divisible by $ x_{{j}1} $. Thus  $ \{s\} \in \mathcal{E} (H) $ with label $ \{x_{{j}t}~|~1 \leq t \leq {\max\{k,{w_{j}}\}} \}. $  
\vspace*{0.2cm}\\	 
Now, we assume  $x_j$ is not a leaf.    Since $x_j$ is not a source,  by our assumption,  either $w_j \neq k$ or $x_i$ is a leaf. If $x_i$ is a leaf, then by the above argument, $m_s$ is the only minimal generator of $I(D)^{\mathcal{P}}$ which is divisible by $x_{i1}$  and therefore $ \{s\} \in \mathcal{E} (H)  $ with label $ \{x_{{i}t}~|~1 \leq t \leq {\min\{k,{w_{i}}\}} \}. $  
If $w_j \neq k$, then ${min\{k,{w_{j}}\}}   < {max\{k,{w_{j}}\}}$.
Here  $\displaystyle{ m_s= \prod_{l=1}^{min\{k,{w_{i}}\}}x_{{i}l}\prod_{t=1}^{max\{k,{w_{j}}\}}x_{{j}t}} $.  Since $|N_D^-(x_j)| = 1  $,    $ m_s $
 is the only minimal generator of $ I(D)^{\mathcal{P}} $
 which is divisible by $ x_{j({\min\{k,{w_{j}}\}+1})} $. Thus $ \{s\} \in \mathcal{E} (H)  $ with label $$ \{x_{{j}t}~|~{\min\{k,{w_{j}}\}+1} \leq t \leq {\max\{k,{w_{j}}\}} \}. $$ Therefore, $\{s\}  \in   E(H(D))$ for all $1 \leq s \leq r$. Hence, by Proposition \ref{labeled}, we get
 \vspace{0.2cm}\\
$\displaystyle{  \reg(R/I(D)) = |X(H)| - |V (H)| = \sum_{\substack{x_i \in V(D)\\~\mbox{and}~N^{-}_D(x_i) \neq \emptyset}} {{\max\{k,{w_{i}}\}}} + \sum_{\substack{x_i \in V(D)\\~\mbox{and}~N^{-}_D(x_i) = \emptyset}} {{{w_{i}}}}-|E(D)|}  $ 
 \vspace{0.2cm}\\	
$\hspace*{0cm}\mbox{and}~\displaystyle{  \pd(R/I(D)) =    |V (H)|  = |E(D)|}  .$
\end{proof}

\begin{remark}
	As an application of Theorem \ref{atmost-one},  we  can get regularity and projective dimension of a large class of $k$-edge weighted VEW graphs such as naturally oriented paths, naturally oriented cycles, rooted forests, naturally oriented unicyclic graphs with a unique cycle, etc.
\end{remark}
In the remainder of this section, we study the regularity of $k$-edge weighted VEW graphs when all the vertices of $ V^* $ are sinks. 
\begin{remark}\label{equalweight}
Let $ D $ be a $k$-edge weighted VEW graph.  If $x_i \in V^*$ is a sink vertex, then by Remark \ref{rem-sink-source}, $w(x_i) > k$ and  $ \max\{k,w(x_i)\} = w(x_i)$. If $x_i \notin V^*$, then $w(x_i)=k$ or $x_i$ is a source vertex.
In case of $w(x_i)=k$, we have $\min\{k,w(x_i)\} = \max\{k,w(x_i)\} = w(x_i) = k$ and if $x_i$ is a source vertex, then by Remark \ref{rem-sink-source}, $\min\{k,w(x_i)\} = w(x_i)$.
\end{remark}
\hspace*{0.5cm}Now, we consider $ D $ to be a $k$-edge weighted VEW graph such that all the vertices of $ V^* $ are sinks. By Remark \ref{equalweight}, for every $x_i \in V(D)$,  the exponent of $x_i = w(x_i)$ in all the minimal generators of $I(D)$ that involves $x_i$. Therefore
$I(D)  =  (x_p^{w(x_p)} x_q^{w(x_q)}~|~(x_p,x_q) \in E(D) ).$  This class of ideal is also mentioned in \cite[Question 5.20]{kara}.

\begin{example}
Consider the  $5$-edge weighted VEW path $D_1$ and the $8$-edge weighted VEW odd cycle $D_2$  as in Figure \ref{fig.1}.  Here $V^*(D_1)= \{x_2  \}$ and $ V^*(D_2)= \{x_1,x_3\}$. Note that $x_2$ is sink in $D_1$, and  $x_1$ and $x_3$ are sinks in $D_2$.  
Therefore, $I(D_1) =(x_1^5x_2^9,x_3^2x_{2}^9,x_3^2x_4^5)$ 
 and  $I(D_2) =(x_2^8x_1^9,x_2^8x_{3}^{10},x_4^8x_3^{10},x_5^5x_4^8,x_5^5x_1^9)$.
 
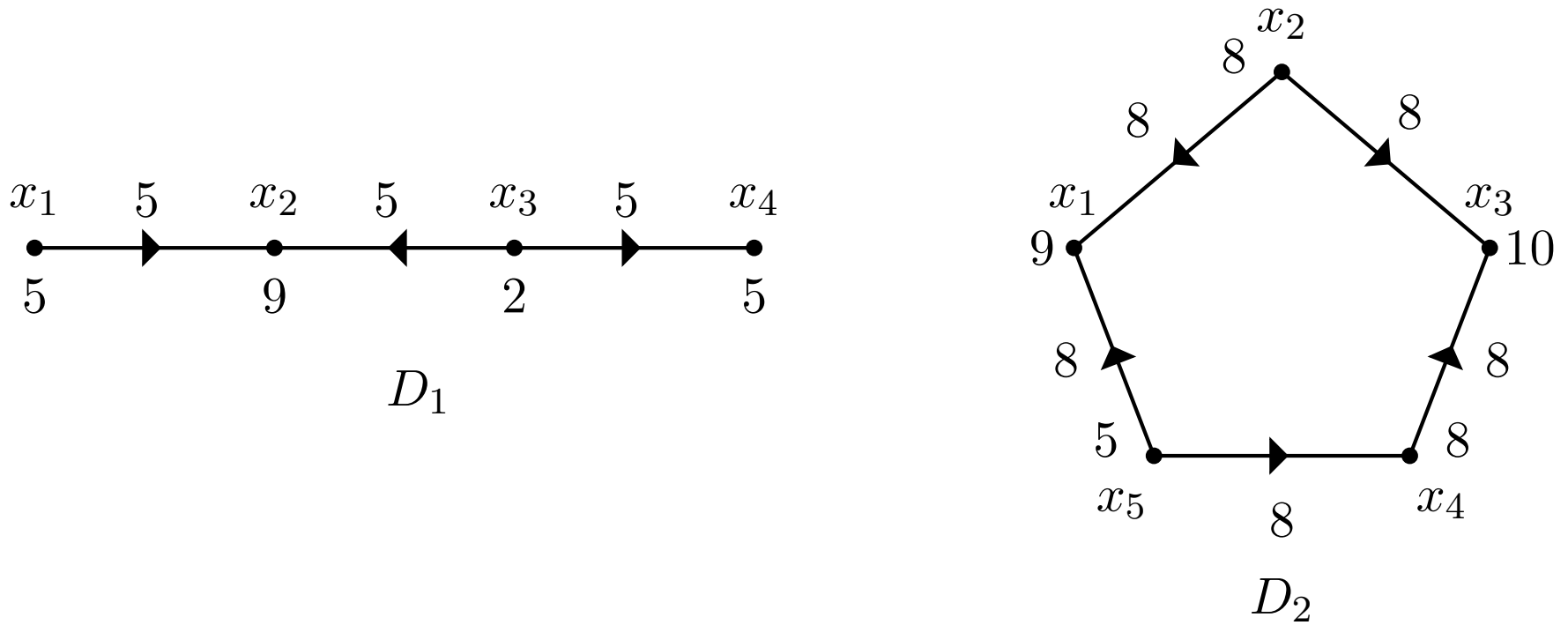
\begin{figure}[!ht]  
	\begin{tikzpicture}[scale=1.2]
	\begin{scope}[thick, every node/.style={sloped,allow upside down}] 
	\definecolor{ultramarine}{rgb}{0.07, 0.04, 0.56} 
	\definecolor{zaffre}{rgb}{0.0, 0.08, 0.66}

\node at (2.8,-0.3) {$x_5$};
\node at (4.8,-0.3) {$x_4$};
\node at (2.5,1.6) {$x_{1}$};
\node at (5.1,1.6) {$x_3$};
\node at (3.8,2.7) {$x_2$};

\node at (2.7,0.1) {$5$};  
\node at (4.9,0.1) {$8$};
\node at (2.3,1.3) {$9$};
\node at (5.35,1.3) {$10$};
\node at (3.5,2.5) {$8$}; 

\node at (2.9,2.1) {$8$};  
\node at (4.6,2.15) {$8$};
\node at (5.15,0.6) {$8$};
\node at (3.8,-0.4) {$8$};
\node at (2.45,0.6) {$8$};

\draw[fill] (3,0) --node {\midarrow}(4.6,0);
\draw[fill] (4.6,0) --node {\midarrow}(5.1,1.3);
\draw[fill] (3,0) --node {\midarrow}(2.5,1.3);
\draw[fill] (3.8,2.4) --node {\midarrow}(5.1,1.3);
\draw[fill] (3.8,2.4) --node {\midarrow}(2.5,1.3);

  \draw [fill] [fill] (3,0) circle [radius=0.04];
  \draw [fill] [fill] (4.6,0) circle [radius=0.04];
  \draw [fill] [fill] (5.1,1.3) circle [radius=0.04];
  \draw [fill] [fill] (3.8,2.4) circle [radius=0.04];
  \draw [fill] [fill] (2.5,1.3) circle [radius=0.04];

\node at (3.8,-0.9) {$D_2$};

\node at (0.5,1.6) {$x_4$};
\node at (-1,1.6) {$x_3$};
\node at (-2.5,1.6) {$x_2$};
\node at (-4,1.6) {$x_1$};

\node at (0.5,1) {$5$};
\node at (-1,1) {$2$};
\node at (-2.5,1) {$9$};
\node at (-4,1) {$5$};

\node at (-1.8,1.6) {$5$};
\node at (-0.3,1.6) {$5$};
\node at (-3.3,1.6) {$5$};

  \draw [fill] [fill] (-2.5,1.3) circle [radius=0.04];
  \draw [fill] [fill] (-4,1.3) circle [radius=0.04];
  \draw [fill] [fill] (-1,1.3) circle [radius=0.04];
  \draw [fill] [fill] (0.5,1.3) circle [radius=0.04];

\draw[fill] (-4,1.3) --node {\midarrow}(-2.5,1.3);
\draw[fill] (-1,1.3) --node {\midarrow}(-2.5,1.3);
\draw[fill] (-1,1.3) --node {\midarrow}(0.5,1.3);

\node at (-1.6,0.4) {$D_1$};
	\end{scope}
	\end{tikzpicture}
		\caption{Two  $k$-edge weighted VEW graphs}\label{fig.1} 
\end{figure}
\end{example}

Next, we relate the regularity of the edge ideal of $k$-edge weighted VEW path or cycle with the regularity of the edge ideal of its underlying graph when the vertices of $ V^*$ are sinks. First, we compute the regularity of the edge ideals of the $k$-edge weighted VEW paths when the vertices of $V^*$ are sinks.

\begin{theorem} \label{ClassII}
	Let $ D $ be a $k$-edge weighted VEW path of length $ n  $  with underlying graph $ G= P_n=x_0 x_1 \ldots x_n $ and such that the vertices of $ V^* $ are sinks. Let $ w_i=w(x_i)  $ for $x_i \in V(D).$ Then

	\begin{equation*}
		  \reg(I(D))= \hspace{.3cm}
		\begin{cases}
			&\reg(I(G)) + \displaystyle{  \displaystyle{ \sum_{x_i \in V(D)}}(w_i - 1) } ~\mbox{  when} ~  n\equiv 1,2~(\mbox{mod}~3),\\
			&\reg(I(G)) + \displaystyle{  \displaystyle{ \sum_{x_i \in V(D)\setminus\{x_j\}}}(w_i - 1) } ~\mbox{when}~ n\equiv 0~(\mbox{mod}~3),
		\end{cases}
	\end{equation*}
\hspace*{4cm}	where $x_j$ is such that $w_j = \min\{w_0,w_3,\ldots,w_{n-3},w_{n}\}$.
\end{theorem}

\begin{proof}
If $A_D=  \reg(I(G)) + \displaystyle{  \displaystyle{ \sum_{x_i \in V(D)}}(w(x_i) - 1) }  $ for $  n=|E(D)|\equiv 1,2~(\mbox{mod}~3)  $ and $B_D=\reg(I(G)) + \displaystyle{  \displaystyle{ \sum_{x_i \in V(D)\setminus\{x_j\}}}(w(x_i) - 1) }$, where  $w(x_j) = \min\{w_0,w_3,\ldots,w_{n-3},w_{n}\}$, for $  n=|E(D)|\equiv 0~(\mbox{mod}~3)  $, then it is sufficient  to prove that
 $$\reg(I(D)) =   \frac{1}{3}[2 -\alpha^{n} -(\alpha^2)^{n}]A_D  + \frac{1}{3}[1 +\alpha^n + (\alpha^2)^n]B_D,$$ where $\alpha$ is the cube root of unity, $\alpha \neq 1$, $1+\alpha + \alpha^2 = 0$ and $\alpha^3 = 1$.
Note that  $\reg(I(D)) = A_D  $  when  $  n\equiv 1,2~(\mbox{mod}~3)  $  and $\reg(I(D)) = B_D  $  when  $  n\equiv 0~(\mbox{mod}~3)  $. Now, we argue by induction on $n = |E(D)|$. If $n=|E(D)| =  1$, then $I(D) = (x_0^{w_0}x_1^{w_1})$ and $\reg(I(D)) = w_0+w_1 = 2 + (w_0-1) + (w_1-1)  =  \reg(I(G)) + (w_0-1) + (w_1-1)=A_D$. Assume that the hypothesis is true for $n=|E(D)| \leq m$ for some $m > 1$. Then we want to show that the hypothesis is true for $n=|E(D)| = m+1.$ Now, we will take three possible cases.       
 
\textbf{Case (1)} Suppose $n=|E(D)| = m+1   \equiv 1~(\mbox{mod}~3)$. Then $ n=m+1 = 3t+1$ for some $t \geq 1$. Let $ D_1 = D\setminus \{x_{3t+1}\} $, $ D_2 = D\setminus \{ x_{3t-1},x_{3t},x_{3t+1} \} ,$ $ G_1 = G\setminus \{x_{3t+1}\} $ and $ G_2 = G\setminus \{ x_{3t-1},x_{3t},x_{3t+1} \} ,$ i.e., $ G_1 $ and $ G_2 $ are the corresponding underlying graphs of $D_1 $ and  $ D_2 $, respectively.
Then
$I(D)=(x_0^{w_{0}}x_1^{w_{1}},x_1^{w_{1}}x_2^{w_{2}},\ldots,x_{3t}^{w_{3t}}x_{3t+1}^{{w_{3t+1}}}).$  Let $ J = (x_{3t}^{w_{3t}}x_{3t+1}^{{w_{3t+1}}}) $ and $ K=I(D_1). $ Since $ J $ has linear resolution, then by Theorem \ref{betti.1}, $ I(D) = J + K $ is a Betti splitting and $ \reg(J)= w_{3t} + {w_{3t+1}} $. As $ D_1 $ is the $k$-edge weighted VEW path of length $ 3t \equiv 0~(\mbox{mod}~3)$, by induction hypothesis, we have $$ \reg(K)=\reg(I(D_1))= B_{D_1} =  \reg(I(G_1)) +   \displaystyle{ \sum_{x_i \in V(D_1)\setminus \{x_j\} }}(w_i - 1),$$ where  $w(x_j) = \min\{w_0,w_3,\ldots,w_{3t}\}$. Note that $ J \cap K= JL $ with $ L= (I(D_2),x_{3t-1}^{w_{3t-1}}) .$ Since $ D_2 $ is the $k$-edge weighted VEW path of length $  3t-2 \equiv 1~(\mbox{mod}~3)$, by induction hypothesis, we have $\reg(I(D_2))= A_{D_2}= \reg(I(G_2)) +  \displaystyle{ \sum_{x_i \in V(D_2)} (w_i - 1) }  .$  By Proposition \ref{path}, $ \reg(I(G))= \reg(I(G_2)) + 1 .$ Thus by Lemma \ref{reg.2}, we have   
\begin{align*}
	\reg(J\cap K) - 1&= \reg(J) + \reg(L)  - 1\\
	&=\reg(J) + \reg(I(D_2)) + w_{3t-1} - 2\\
	&= w_{3t} + {w_{3t+1}} + \reg(I(G_2)) +   \displaystyle{ \sum_{x_i \in V(D_2)} (w_i - 1) } + w_{3t-1} - 2\\
	&= \reg(I(G_2)) + 1 + \displaystyle{  \displaystyle{ \sum_{x_i \in V(D) }}(w_i - 1) }\\
	&=  \reg(I(G)) + \displaystyle{  \displaystyle{ \sum_{x_i \in V(D) }}(w_i - 1) }.
\end{align*} 
By Lemma \ref{H<G}, $\reg(I(G_1)) \leq \reg(I(G))	.$  Then by Corollary \ref{betti.2}, we get    
\begin{align*}
	\reg(I(D)) &= \max\{ \reg(J),~\reg(K), ~\reg(J\cap K) - 1\}\\&= \max\{ w_n + 1,~\reg(I(G_1)) +   \displaystyle{ \sum_{x_i \in V(D_1) }}(w_i - 1) ,~ \reg(I(G)) + \displaystyle{ \sum_{x_i \in V(D)}}(w_i - 1)  \} \\
	&=\reg(I(G)) + \displaystyle{ \sum_{x_i \in V(D)}}(w_i - 1)\\
    &=A_D.
\end{align*}

\textbf{Case (2)} Suppose $n=|E(D)| = m+1   \equiv 2~(\mbox{mod}~3)$. Then $n=m+1 = 3t+2$ for some $t \geq 0$. By a similar argument as in Case (1), we can show that $\reg(I(D))  =  A_D  $.    
 
\textbf{Case (3)} Supoose $n=|E(D)| = m+1  \equiv 0~(\mbox{mod}~3)$. Then $n=m+1 = 3t+3$ for some $t \geq 0$. Let $ D_1 = D\setminus \{x_{3t+3}\} $, $ D_2 = D\setminus \{ x_{3t+1},x_{3t+2},x_{3t+3} \} ,$ $ G_1 = G\setminus \{x_{3t+3}\} $ and $ G_2 = G\setminus \{ x_{3t+1},x_{3t+2},x_{3t+3} \} ,$ i.e., $ G_1 $ and $ G_2 $ are the corresponding underlying graphs of $D_1 $ and  $ D_2 $, respectively. Then
 $I(D)=(x_0^{w_{0}}x_1^{w_{1}},x_1^{w_{1}}x_2^{w_{2}},\ldots,x_{3t+2}^{w_{3t+2}}x_{3t+3}^{{w_{3t+3}}}).$  Let $ J = (x_{3t+2}^{w_{3t+2}}x_{3t+3}^{{w_{3t+3}}}) $ and $ K=I(D_1). $ Since $ J $ has linear resolution, $ I(D) = J + K $ is a Betti splitting and $ \reg(J)= w_{3t+2} + {w_{3t+3}} $. As $ D_1 $ is the $k$-edge weighted VEW path of length $ 3t+2 \equiv 1~(\mbox{mod}~3)$, by induction hypothesis, we have $$ \reg(K)=\reg(I(D_1))=A_{D_1}= \reg(I(G_1)) +   \displaystyle{ \sum_{x_i \in V(D_1) }}(w_i - 1).$$ Note that $ J \cap K= JL $, where $ L= (I(D_2),x_{3t+1}^{3t+1}) .$ Since $ D_2 $ is the $k$-edge weighted VEW path of length $ 3t \equiv 0~(\mbox{mod}~3) $, by induction hypothesis, we have $$\reg(I(D_2))=B_{D_2}= \reg(I(G_2)) +  \displaystyle{ \sum_{x_i \in V(D_2)\setminus \{x_p\}} (w_i - 1) },$$ where  $w_p = \min\{w_0,w_3,\ldots,w_{3t}\}$. By Proposition \ref{path}, $ \reg(I(G))= \reg(I(G_2)) + 1 .$ Thus, by Lemma \ref{reg.2}, we have   
 \begin{align*}
 	\reg(J\cap K) - 1&= \reg(J) + \reg(L)-1\\
 	&=\reg(J) + \reg(I(D_2)) + w_{3t+1} - 2\\
 	&= w_{3t+2} + {w_{3t+3}} + \reg(I(G_2)) +   \displaystyle{ \sum_{x_i \in V(D_2)\setminus \{x_p\}} (w_i - 1) } + w_{3t+1} - 2\\
 	&= \reg(I(G_2)) + 1 + \displaystyle{  \displaystyle{ \sum_{x_i \in V(D)\setminus \{x_p\} }}(w_i - 1) }\\
 	&=  \reg(I(G)) + \displaystyle{  \displaystyle{ \sum_{x_i \in V(D)\setminus \{x_p\} }}(w_i - 1) }.
 \end{align*} 
 By Proposition \ref{path}, $\reg(I(G_1)) = \reg(I(G))	.$  Then by Corollary \ref{betti.2}, we get    
 \begin{align*}
 	\reg(I(D)) &= \max\{ \reg(J),~\reg(K), ~\reg(J\cap K) - 1\}\\&= \max\{ w_{3t+2} + {w_{3t+3}},~\reg(I(G)) +   \displaystyle{ \sum_{x_i \in V(D)\setminus \{x_{3t+3}\} }}(w_i - 1) ,\\
  &\hspace*{0.3cm}  \reg(I(G)) + \displaystyle{  \displaystyle{ \sum_{x_i \in V(D)\setminus \{x_p\} }}(w_i - 1) }  \} \\
 	&=\reg(I(G)) + \displaystyle{  \displaystyle{ \sum_{x_i \in V(D)\setminus\{x_j\}}}(w_i - 1) },
    \end{align*}
    where $w(x_j)=\min\{w_{3t+3},w_p\}  = \min\{w_0,w_3,\ldots,w_{3t}, w_{3t+3}\}$. Hence, $\reg(I(D))=B_D$.

\end{proof}

Next, we compute the regularity of the edge ideals of $k$-edge weighted VEW cycles when the vertices of $V^*$ are sinks.

\begin{theorem} \label{cycle1}
	Let $ D $ be a  $k$-edge weighted VEW cycle of length $ n $ for $ n\equiv 0,1~(\mbox{mod}~3)$  with  underlying graph  $ G= C_n=x_1\ldots x_nx_1$ and such that the vertices of $ V^* $ are sinks. Then $$\reg(I(D)) = \reg(I(G)) + \displaystyle{ \sum_{x_i \in V(D)}}(w(x_i) - 1).$$   
\end{theorem}

\begin{proof}
Let $ V(D)=\{x_1,\ldots,x_n\} $ and let $ w_i=w(x_i)  $ for $i=1, \ldots, n.$ Since vertices of $ V^* $ are sinks, we have $I(D) = (x_i^{w_i}   x_j^{w_j}  ~|~  (x_i,x_j) \in E(D))$. Without loss of generality, let  $ x_p  \neq x_1,x_n$.  
	Let $ D_1 = D \setminus \{x_p\} ,$ $ D_2 = D \setminus N_D[x_p] ,$ $ G_1=G \setminus \{x_p\} $ and  $ G_2 = G \setminus N_G[x_p] ,$ i.e., $ G_1 $ and $ G_2 $ are the corresponding underlying graphs of $D_1 $ and  $ D_2 $, respectively.  
	By Proposition \ref{path}, $\reg(I(G))   =  \reg(I(G_1))=\reg(I(G_2)) + 1    $ except $n=3,4$ 
	and $\reg(I(G))   =  \reg(I(G_1))=2$ for $n=3,4$. Now, consider the exact sequence
	\begin{equation}\label{2}
		0 \longrightarrow {\frac{R}{I(D):x_p^{w_p} }}(-w_p) \xrightarrow{.x_p^{w_p}} {\frac{R}{I(D) }} \longrightarrow {\frac{R}{(I(D),x_p^{w_p} )}} \longrightarrow 0.
	\end{equation}		
	Here $(I(D),x_p^{w_p})= (I(D_1),x_p^{w_p})$, where $ D_1 = D \setminus \{x_p\} $ is a $k$-edge weighted VEW path of length $ n-2 $ with end vertices $x_{p-1}$ and $x_{p+1}$. Note that $ n-2\nequiv 0~(\mbox{mod}~3)$.  Thus by Lemma \ref{reg.2}(1) and Theorem \ref{ClassII}, we have
	\begin{align*}
		\reg((I(D),x_p^{w_p}))&=\reg(I(D_1))+ {w_p} - 1\\
		&=\reg(I(G_1)) +   \displaystyle{ \sum_{x_i \in V(D_1) }}(w_i - 1)+ {w_p} - 1\\
        &= \reg(I(G)) +   \displaystyle{ \sum_{x_i \in V(D)}}(w_i - 1).
	\end{align*}
	Here $I(D):x_p^{w_p}=(I(D_2),x_{p-1}^{w_{p-1}},x_{p+1}^{w_{p+1}})$ for $n \neq 3,4$ and $I(D):x_p^{w_p}=
	(x_{p-1}^{w_{p-1}},x_{p+1}^{w_{p+1}})$ for $n=3,4.$ For $n  \neq 3, 4$, $ D_2  $ is a $k$-edge weighted VEW path of length $ n-4 $. By Lemma \ref{reg.2}(1) and Theorem    \ref{ClassII}, we have
	\begin{align*}
		\reg((I(D):x_p^{w_p})(-{w_p}))&=\reg(I(D_2)) +  w_{p-1} + w_{p+1} - 2 +  {w_p}\\
		&\leq \reg(I(G_2)) +   \displaystyle{ \sum_{x_i \in V(D_2) }}(w_i - 1)+  w_{p-1} + w_{p+1} - 2 +  {w_p}\\
		&=\reg(I(G_2)) + 1 +   \displaystyle{ \sum_{x_i \in V(D)}}(w_i - 1)\\
		&=\reg(I(G)) +   \displaystyle{ \sum_{x_i \in V(D)}}(w_i - 1).
	\end{align*}
	
	For $n=3,4,$ $\reg((I(D):x_p^{w_p})(-{w_p}))= w_{p-1} + w_{p+1} - 1 + {w_p}\\
	\hspace*{6.1cm}= 2 + (w_{p-1}-1) + (w_{p+1} - 1) + ({w_p}-1)\\
	\hspace*{6.1cm} \leq   \reg(I(G)) +    \displaystyle{ \sum_{x_i \in V(D) }}(w_i - 1).$ 
	
	By Lemma \ref{reg.4} and the exact sequence (\ref{2}), we get
	\hspace*{4cm} $$ \reg(I(D)) \leq \max\{\reg((I(D):x_p^{w_p})(-{w_p})),\reg((I(D),x_p^{w_p} ))\} .$$
	Since $\reg((I(D):x_p^{w_p})(-{w_p})) - 1 \neq \reg((I(D),x_p^{w_p})),$ by Lemma  \ref{reg.4} and the exact \\sequence (\ref{2}), we have
	\begin{align*}
		\reg(I(D)) &=\reg((I(D),x_p^{w_p})) = \reg(I(G)) +   \displaystyle{ \sum_{x_i \in V(D)}}(w_i - 1).
	\end{align*}
    
\end{proof}

\begin{theorem} \label{cycle2}
	Let $ D $ be a $k$-edge weighted VEW cycle of length $ n $ for $ n\equiv 2~(\mbox{mod}~3)$  with  underlying graph  $ G= C_n=x_1\ldots x_nx_1$ and such that the vertices of $ V^* $ are sinks. Then $$\reg(I(D)) = \reg(I(G)) +   \displaystyle{ \sum_{x_i \in V(D)}}(w(x_i) - 1).$$  
\end{theorem}

\begin{proof}
Let $ V(D)=\{x_1,\ldots,x_n\} $ and let $ w_i=w(x_i)  $ for $i=1, \ldots, n.$ We argue by induction on the number $m$ of vertices of $V^*$.  If $m=0$, then the result follows trivially. Now, let $m \geq 1$ and  without loss of generality, that $x_p \neq x_1, x_n$ is in $V^*$. Note that $w_p >k$.
	Let $ D_1 = D \setminus \{x_p\} $ and  $ G_1=G \setminus \{x_p\} ,$      i.e., $ G_1 $ is  the corresponding underlying graph of $D_1 $.      
	By Proposition \ref{path}, we have $ \reg(I(G)) = \reg(I(G_1)) + 1.$ 	
	Now, consider the exact sequence  
	\begin{equation}\label{3}
		0 \longrightarrow {\frac{R}{I(D):x_p^{w_p - k} }}(-(w_p-k)) \xrightarrow{.x_p^{w_p - k}} {\frac{R}{I(D) }} \longrightarrow {\frac{R}{(I(D),x_p^{w_p - k}) }} \longrightarrow 0 .
	\end{equation}
	Here $(I(D),x_p^{w_p - k})= (I(D_1),x_p^{w_p - k})$, where  $ D_1 = D \setminus \{x_p\} $ is a $k$-edge weighted VEW path of length $ n-2 $ with end vertices $x_{p-1}$ and $x_{p+1}$. Note that $n-2\equiv 0~(\mbox{mod}~3)$. Then by Theorem \ref{ClassII}, we have $   \reg(I(D_1))=  \reg(I(G_1)) +   \displaystyle{ \sum_{x_i \in V(D_1)\setminus \{x_j\}}}(w_i - 1)$,
	 where  $ w(x_j) = \min\{\ldots, w(x_{p-4 \mod~n}),  w(x_{p-1}),     w(x_{p+1}), w(x_{p+4 \mod~n}),\ldots\}.$    
	Thus by  Lemma \ref{reg.2}(1), we have  
	\begin{align*}
		\reg((I(D),x_p^{w_p - k}))&= \reg(I(D_1))+ (w_p -k) -1\\
		&=  \reg(I(G_1)) +   \displaystyle{ \sum_{x_i \in V(D_1) \setminus \{x_j\}}}(w_i - 1)+ (w_p -k) -1\\  
		&=  \reg(I(G)) +   \displaystyle{ \sum_{x_i \in V(D)\setminus \{x_j\}}}(w_i - 1) - k - 1.
	\end{align*}
	
	Here $I(D):x_p^{w_p - k}=I(D_3)$, where $ D_3 $ is a $k$-edge weighted VEW cycle with $m-1$  vertices in $V^*(D_3)$. Notice that $w(x_p)=k$ in $D_3$ and $w(x_p) > k$ in $D$. By using the induction hypothesis,$$ \reg(I(D_3))=  \reg(I(G)) +   \displaystyle{ \sum_{x_i \in V(D_3)}}(w_i - 1)=  \reg(I(G)) +   \displaystyle{ \sum_{x_i \in V(D)\setminus \{x_p\}}}(w_i - 1)+k-1.$$
	Then   $\reg((I(D):x_p^{w_p - k})(-(w_p-k)))=\reg(I(D_3)) + w_p - k =\reg(I(G)) +   \displaystyle{ \sum_{x_i \in V(D)}}(w_i - 1).$ Thus, by applying Lemma \ref{reg.4} on exact sequence (\ref{3}), we have
	$$ \reg(I(D)) \leq \max\{\reg((I(D):x_p^{w_p - k})(-(w_p-k))),\reg((I(D),x_p^{w_p - k}))\} .$$  
	Since $\reg((I(D):x_p^{w_p - k})(-(w_p-k))) - 1  \neq  \reg((I(D),x_p^{w_p - k})),$ by applying Lemma \ref{reg.4} on exact sequence (\ref{3}), we get 
	$$ \reg(I(D))= \reg((I(D):x_p^{w_p - k})(-(w_p-k)))=\reg(I(G)) +   \displaystyle{ \sum_{x_i \in V(D)}}(w_i - 1) . $$  	 
	
\end{proof}

\section{Discussion}\label{sec7}
We have studied the regularity of $k$-edge weighted VEW graphs with at most one edge oriented into each of its vertex and  $k$-edge weighted VEW paths or cycles assuming non-source vertices with weight $\neq k  $  are sinks. There are different classes of $k$-edge weighted VEW graphs, whose regularity need to be studied.
For weighted oriented graphs, i.e. $1$-edge weighted VEW graphs, symbolic powers have been studied in \cite{kanoy1, mandal2,mandal3,mandal4}. Now for arbitrary $k$, the following question arises naturally.

\begin{question}
\begin{enumerate}
    \item[(a)] Compare the ordinary and symbolic powers of edge ideals of different classes of $k$-edge weighted VEW graphs.  

    \item[(b)] Compute the regularity of powers of edge ideals of different classes of $k$-edge weighted VEW graphs.              

\end{enumerate} 
\end{question}

\begin{example}
Consider the  VEW graph $D$ with different weights on the edges as in Figure \ref{fig.4}.  Then  $I(D) =(x_1^3x_2^6,x_2^2x_{3}^3,x_4^2x_3^5,x_5^2x_4^7,x_5^2x_1^3)$.  

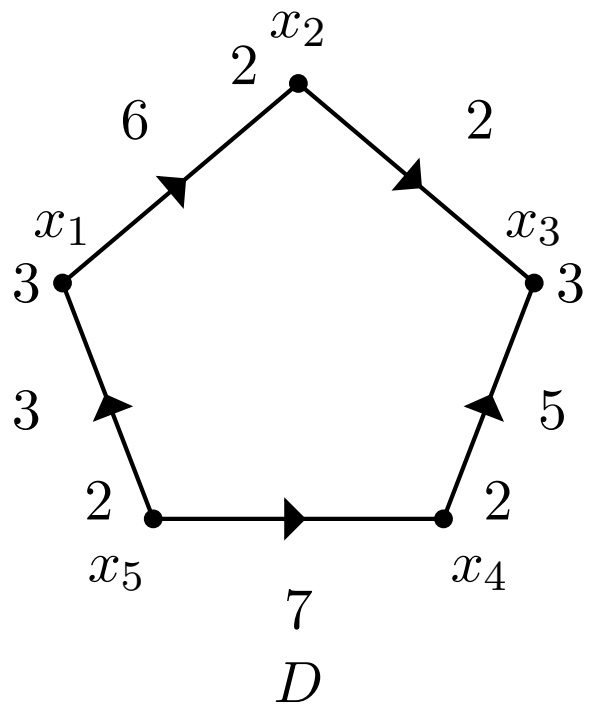
\begin{figure}[!ht]  
	\begin{tikzpicture}[scale=1.2]
	\begin{scope}[thick, every node/.style={sloped,allow upside down}] 
	\definecolor{ultramarine}{rgb}{0.07, 0.04, 0.56} 
	\definecolor{zaffre}{rgb}{0.0, 0.08, 0.66}

\node at (2.8,-0.3) {$x_5$};
\node at (4.8,-0.3) {$x_4$};
\node at (2.5,1.6) {$x_{1}$};
\node at (5.1,1.6) {$x_3$};
\node at (3.8,2.7) {$x_2$};

\node at (2.7,0.1) {$2$};  
\node at (4.9,0.1) {$2$};
\node at (2.3,1.3) {$3$};
\node at (5.3,1.3) {$3$};
\node at (3.5,2.5) {$2$}; 

\node at (2.9,2.2) {$6$};  
\node at (4.8,2.2) {$2$};
\node at (5.2,0.6) {$5$};
\node at (3.8,-0.5) {$7$};
\node at (2.3,0.6) {$3$};

\draw[fill] (3,0) --node {\midarrow}(4.6,0);
\draw[fill] (4.6,0) --node {\midarrow}(5.1,1.3);
\draw[fill] (3,0) --node {\midarrow}(2.5,1.3);
\draw[fill] (3.8,2.4) --node {\midarrow}(5.1,1.3);
\draw[fill] (2.5,1.3) --node {\midarrow}(3.8,2.4);

  \draw [fill] [fill] (3,0) circle [radius=0.04];
  \draw [fill] [fill] (4.6,0) circle [radius=0.04];
  \draw [fill] [fill] (5.1,1.3) circle [radius=0.04];
  \draw [fill] [fill] (3.8,2.4) circle [radius=0.04];
  \draw [fill] [fill] (2.5,1.3) circle [radius=0.04];

	\node at (3.8,-0.9) {$D$};

	\end{scope}
	\end{tikzpicture}
		\caption{A VEW cycle of length 5}\label{fig.4} 
\end{figure}
Using Macaulay2 \cite{gs}, the irredundant irreducible decomposition of $I(D)$ is given as
	\begin{align*}
	I(D) 
&=  (x_1^3,x_2^2,x_4^2) \cap 
 (x_1^3,x_3^3,x_4^7) \cap  (x_1^3,x_3^3,x_5^2) \cap  (x_2^6,x_2^2x_3^3, x_3^5,x_5^2)\\
&\hspace{0.4cm}  \cap  (x_2^2,x_4^2,x_5^2)  \cap  (x_1^3, x_2^2, x_3^5, x_4^7)   \\  
&=  (x_1^3,x_2^2,x_4^2) \cap 
 (x_1^3,x_3^3,x_4^7) \cap  (x_1^3,x_3^3,x_5^2) \cap  (x_2^2, x_3^5,x_5^2)\cap  (x_2^6,x_3^3,x_5^2)\\
&\hspace{0.4cm}  \cap  (x_2^2,x_4^2,x_5^2)  \cap  (x_1^3, x_2^2, x_3^5, x_4^7).
	\end{align*}
Observe that $\rad (x_2^2, x_3^5,x_5^2) = (x_2, x_3,x_5)$
and $\rad (x_2^6,x_3^3,x_5^2) = (x_2, x_3,x_5)$. That means the radicals of the elements of the irredundant irreducible decomposition of $I(D)$ are not different.          
\end{example}

In the irredundant irreducible decomposition of edge ideals of  $k$-edge weighted VEW graphs, we have seen that radicals of irreducible ideals are distinct.    
Therefore, for the arbitrary weights on the edges of a VEW graph, we have the following question.         

\begin{question}
Find the  irredundant  irreducible decomposition of $I(D)$, where $D$ is a VEW graph.

\end{question}

{\bf Acknowledgement:} 
Dipak Kumar Pradhan wants to thank Prof. Chanchal Kumar for posing a question about Vertex-edge weighted oriented graph.

\end{document}